\documentclass[10
 pt,twoside]{amsart}
  \usepackage{amsmath,amssymb,amsthm,amscd,latexsym,graphicx,tikz}
 \usepackage{epigraph}
 \usepackage[OT2,OT1]{fontenc}
  \usepackage{enumerate}

\usepackage{graphicx}
\usepackage[all]{xy}

  \entrymodifiers={!!<0pt,0.7ex>+}

\newcommand{\fonction}[5]{\begin{array}{cccc}
#1: & #2 & \longrightarrow & #3 \\
 & #4 & \longmapsto & #5 \end{array}}
 
\newcommand{\fonctionnoname}[4]{\begin{array}{ccc}
 #1 & \longrightarrow & #2 \\
 #3 & \longmapsto & #4 \end{array}}

 \newcommand{\Fl}{\mathbf{Fl} }

 \renewcommand{\U}{\mathbf{U} }

 \newcommand{\obs}{\mathrm{Obs}}

       \newcommand{\B}{\mathbf{B}}

 \renewcommand{\to}{\longrightarrow}

 \newcommand{\A}{\mathbb A}
 \newcommand{\T}{\mathcal T}

 \newcommand{\F}{\mathbb{F}}
 \newcommand{\G}{\mathbf{G}}

 \newcommand{\N}{\mathbb{N}}
 \renewcommand{\P}{\mathbb P}
 \newcommand{\Q}{\mathbb{Q}}
 
 \newcommand{\Z}{\mathbb{Z}}

 \newcommand{\Res}{\mathrm{Res}}

 \newcommand{\GL}{\mathbf{GL}}

 \newcommand{\Sch}{\mathbf{Sch}}
 \newcommand{\Sets}{\mathbf{Sets}}
 
 \newcommand{\Det}{\mathrm{Det}}
 \newcommand{\Aff}{\mathbf{Aff}}

 \newcommand{\End}{\mathrm{End}}

 \newcommand{\Ext}{\mathrm{Ext}}

 \newcommand{\Hom}{\mathrm{Hom}}

 \newcommand{\Ver}{\mathrm{Ver}}
 \newcommand{\Frob}{\mathrm{Frob}}

 \newcommand{\RR}{{\mathbf{R}}}
 
 \newcommand{\TT}{\mathbf{T}}

 \renewcommand{\T}{\mathbf{T}}

\newcommand{\Id}{\mathrm{Id}}

\renewcommand{\Im}{\mathrm{Im}}
 
\newcommand{\nat}{\mathrm{nat}}
\newcommand{\Ker}{\mathrm{Ker}}
 
 \newcommand{\Spec}{\mathrm{Spec}}
 
 \newcommand{\AAut}{\mathbf{Aut}}

 \newcommand{\ver} {\mathrm{ver}}

 \renewcommand{\O}{\mathcal{O}}

\newcommand{\W}{\mathbf{W}}

\newcommand{\frob}{\mathrm{frob}}

\newcommand{\EExt}{\mathbf{Ext}}

\newcommand{\perf}{\mathrm{perf}}
\newcommand{\Sym}{\mathrm{Sym}}

 \theoremstyle{plain}
 \newtheorem{thm}{Theorem}[section]
   \newtheorem*{thm*}{Theorem}
 \newtheorem{defi}[thm]{Definition}

 \newtheorem{prop}[thm]{Proposition}

 \newtheorem{lem}[thm]{Lemma}

 \newtheorem*{theorem-non}{Bloch-Kato conjecture, an equivalent formulation}
\newtheorem*{thmA*}{Theorem A}
\newtheorem*{thmB*}{Theorem B}
\newtheorem*{thmC*}{Theorem C}
\newtheorem*{thmD*}{Theorem D}

 \theoremstyle{remark}
 \newtheorem{rem}[thm]{Remark}

 \newtheorem{mot}[thm]{Motto}

 \newtheorem{rems}[thm]{Remarks}
 \newtheorem{ex}[thm]{Example}

 \newtheorem{exo}[thm]{Exercise}

 \newenvironment{dem}{{\bf Proof.}}{\hfill$\square$}
 
\date{\today}

\begin{document}
\title{Smooth profinite groups, II: the Uplifting Pattern}

\author{Mathieu Florence}
\address{Sorbonne Université and Université Paris Cité, CNRS, IMJ-PRG, F-75005 Paris, France.}
\email{mathieu.florence@imj-prg.fr}



\begin{abstract}
    
This text presents a scheme-theoretic enhancement of the theory of smooth profinite groups and cyclotomic pairs, introduced in \cite{DCF1}. To do so, our main technical tools are Hochschild cohomology of affine group schemes and lifting frobenius of vector bundles. The main contribution of this work is the Uplifting Pattern. It is a natural process, to  lift a given equivariant extension of vector bundles, to  its  $\W_2$-counterpart, upon a  `reasonable' combination of base-change and group-change. This is the key ingredient  to prove the Smoothness Theorem, in the paper  \cite{DCF3}.
\end{abstract}
\maketitle

\tableofcontents
\newpage
\section{Introduction.}

\noindent Our task here is to develop a bunch of  tools, that are used to prove the Smoothness Theorem in \cite{DCF3}. Let $p$ be a prime. Recall that this series of three papers, is motivated by the following belief.\\

\noindent (B): The Bloch-Kato-Milnor Conjecture follows from mod $p^2$ Kummer theory.\\

\noindent In the course of investigating (B), we  laid, in the recent work \cite{DCF1}, a tentative axiomatisation of  Kummer theory. In particular, the notion of a $(1,1)$-cyclotomic pair $(G,\Z/p^2(1))$ was introduced. Here $G$ is a profinite group, and up to isomorphism, $\Z/p^{2}(1)$ is simply  a continuous character $$G \to (\Z/p^{2})^\times  \simeq  \F_p^\times \times (\Z/p),$$ playing the role of the usual cyclotomic character mod $p^2$  (see \cite{DCF1} for details).\\ 
\noindent This text can be thought of  as an enhancement of the theory of cyclotomic pairs, to the scheme-theoretic setting. Roughly speaking, (pro)finite groups $G$ are  upgraded to smooth affine group schemes  $\Gamma$ (typically $\GL_N$, as a group scheme over $\Z$), while  $\F_p$-representations of $G$, are upgraded to $\Gamma$-equivariant vector bundles. For the reader to  get a first glimpse of this more technological, yet more `computable' setting, let us mention  the analogue of the assertion\begin{center}
`a representation $G\xrightarrow{\rho_1} \GL_N(\F_p)$ lifts, to a representation $G\xrightarrow{\rho_2} \GL_N(\Z/p^2)$.'   \end{center}

In that upgraded setting,  it becomes \begin{center}
`a $\Gamma$-vector bundle  over a $\Gamma$-scheme $S$ lifts, to a $\Gamma \W_2$-bundle over $S$.'   
\end{center}

\noindent The objective of this paper is to present the Uplifting Pattern. Given a $\Gamma$-equivariant \textit{extension} $E$ of vector bundles over some $\Gamma$-scheme $S$, this pattern is a natural  process for enabling the existence of a lifting of $E$, to an equivariant   extension $E_2$ of $\W_2$-bundles, prescribing the kernel and cokernel of $E_2$.  This is achieved upon applying two compatible operations: a \textit{base-change} $U \to S$, and a change of the acting group $\Gamma' \to \Gamma$ (the latter is an explicit surjection of group schemes, called \textit{group-change}). See the final part of the introduction for details.\\
At this point, an important observation is in order.  In the Uplifting Pattern, it matters to work over a $p$-torsion-free base $S$ (typically, a smooth $\Z$-scheme), rather than over an $\F_p$-scheme. Here is a concrete reason for this. Working `integrally', some otherwise intricate cohomological computations are considerably simplified, by considering the effect of inverting $p$ on $S$. In a sense made precise in \cite{DCF3}, this allows to reduce the cohomological degree of obstructions to lifting- a gain that turns out to be invaluable.
On the one hand, this  phenomenon explains the very short length of \cite{DCF3}, the third paper of this series, where the Uplifting Pattern is crucially applied, to prove the Smoothness Theorem. On the other hand, it also explains the heavier technical aspect of the present paper: achieving a proper definition of these `integral models', requires a  deeper  understanding of Witt vectors, than while dealing with $\F_p$-schemes.   Also, the notion of a WTF data, featuring the exotic ring scheme $\W_2^{[r]}$ over $\Z$, comes in handy to incorporate an `integral' Frobenius twist.\\

\noindent  We comment on  the contents of this paper, assuming familiarity with  \cite{DCF1} and \cite{DCFL}.

\begin{itemize}
 \item{Hochschild cohomology of affine group schemes is an essential and robust tool. Our reference is \cite{J}. Appropriate recollections on this topic are given in section  \ref{SecHoch}. In section \ref{SecHochComp},   we prove a  vanishing result used in  \cite{DCF3}.} 
\item{In section \ref{seceqAG}, we define our (equivariant) scheme-theoretic framework for Witt Modules. Let us point out that equivariance comes here for free, and that truncated Witt vectors $\W_r$ are to be considered  over  $\Z$. In that setting, recollections on `transfer w.r.t. ring schemes' are provided, in the degree of generality that we actually need.}
 \item {Equivariant extensions of $\RR$-Modules, and their link to Hochschild cohomology. Here `equivariant' usually refers to the action of a smooth affine  group scheme. See Section \ref{secext}.} 
     \item {Cyclotomic pairs and  smooth profinite groups, as introduced in \cite{DCF1}. Recollections on this topic are provided in Section \ref{CycloSmooth}.}  \item{The  notion of  a WTF data is introduced in section  \ref{SecWTF}. It features somewhat exotic  two-dimensional ring schemes over $\Z$, to be thought of as `Frobenius twists of the ring scheme  $\W_2$'. A crucial example of such a gadget is given is section  \ref{SecR2}: for each integer $r \geq 0$, one defines a  ring scheme $\RR:=\W_2^{[r]}$, over $\Z$. If $r=0$,  $\RR$ is just $\W_2$.  Its purpose is suggested by the tentative  \begin{mot}
Over an $\F_p$-base, lifting a vector bundle to an  $\RR$-bundle, is equivalent to lifting  its $r$-th (absolute) frobenius twist, to a $\W_2$-bundle.
\end{mot}}  \item {Splitting schemes of extensions of vector bundles are studied in section \ref{SecAff}. They give rise to natural increasing filtrations (on functions and on vector bundles), indexed by  well-ordered sets, called `good filtrations'.} 
   
\item {How to lift equivariant extensions of vector bundles to their $\RR$-counterparts, is the central  aspect of this paper.  This is achieved by  the so-called `Uplifting Pattern' (section \ref{secUP}), that combines a \textit{base-change} and a  \textit{group-change}. The base-change is the  so-called `Uplifting Scheme',  closely related to a   classical  problem: lifting the frobenius of the mod $p$ reduction of a scheme $S$, to an endomorphism of $S$. This is the object of  section \ref{SecLiftFrob}. The group-change features a  \textit{suspension} process, for lifting  Hochschild $1$-cocycles (or equivalently, for trivialising  Hochschild $2$-cocycles)- see Section  \ref{SecSuspension}.}

\end{itemize}

\section{Notation, conventions.}
\noindent The notation $A:=E$ means that $A$ is \textit{defined as} the expression $E$. \\
The letter $p$ denotes a prime number.\\ Oddly enough, the parity of $p$ plays no role. Even better: $p$ can be arbitrary! ;) \\
Unless specified otherwise, rings are unital and commutative. \\Schemes, and morphisms between them, are quasi-compact and quasi-separated. 

\subsection{The notation $*$.}\hfill\\
\noindent  In this text, we use the symbol $*$ to denote an integer, or an object,  whose name it is superfluous to mention. For instance,  \[ 0 \to A \to * \to * \to  B \to 0\]denotes  a $2$-extension of $B$ by $A$ (in an  additive category), whose middle terms need not be specified. 
\subsection{Equivariance.\\}

In \cite{DCF1}, equivariance is always w.r.t.  actions of (pro)finite groups $G$. To give a meaning to the Uplifting Pattern,  these are upgraded here, to (affine and smooth) group schemes $\Gamma$ over an affine base- typically $\Spec(\Z)$. To keep things simple,  pro-algebraic group schemes are avoided. This is possible at the  cost  of  (gladly!) paying a finer attention, to those surjections $\Gamma' \xrightarrow{\gamma} \Gamma$  that would be the transition maps in such projective limits. In many cases, but not always, $\Ker(\gamma)$ is a split unipotent  affine $\Z$-group scheme.\\
Let $$\Gamma \to S$$ denote an affine smooth group scheme, over an affine scheme $S$. To fix ideas, the reader may assume $S=\Spec(\Z)$. Notation-wise, if $\mathcal T$ is a `type' of algebro-geometric objects $X$ over $S$ (e.g. an $S$-scheme, or a Witt vector bundle over an $S$-scheme), we denote by $\Gamma \mathcal T$ the `$\Gamma$-equivariant' version of $\mathcal T$. It consists of objects $X$ of type $\mathcal T$,  with the extra datum of an action of $\Gamma$.  At the level of functors of points, this action amount to the following. For every morphism $S' \to S$, an action of $\Gamma(S')$ is given on the $S'$-structure  $X \times_S S'$, functorially in $S' \to S$.
\subsection{Geometric triviality.\\}
\noindent  Generally speaking, `a geometrically trivial structure' refers to `an equivariant structure, that is trivial upon forgetting equivariance'.

\section{Admissible extensions and Hochschild cohomology.}\label{SecHoch}

 In this section,  $S$  is a scheme. 
\begin{defi}
   An $S$-functor is a contravariant functor \[ \{ \Sch/S\} \to \{\Sets\}.\]  An $S$-group functor (or $S$-group) is a contravariant functor \[ \{ \Sch/S\} \to \{\mathbf{Grp}\}.\] 
\end{defi}

\begin{rem}
    The functor of points of an $S$-scheme (resp. $S$-group scheme) is an $S$-functor (resp. $S$-group).
\end{rem}
 \begin{defi}
 Let $\Gamma$ be an $S$-group.
     A $(\Gamma,S)$-group  is an $S$-group  $H$, equipped with a $\Gamma$-action, by $S$-group automorphisms.
     \end{defi}

  \begin{defi}\label{defi admiss}

   Consider a sequence of $S$-groups, $$\mathcal E: 1 \to A \xrightarrow{i}   B  \xrightarrow{\pi}   C\to 1.$$ Say that $\mathcal E$ is an admissible extension if the three conditions below hold. \begin{enumerate}
       \item{$i$ is injective,}  \item{$(\pi \circ i)=1$,}   \item{There exists a morphism of $S$-functors $\tau: C \to B$,  such that $\pi  \tau=\Id_C.$ }
   \end{enumerate}
      
  \end{defi}

\begin{rem}
    In item (3), it is not required that the section $\tau$ be a homomorphism.
\end{rem}

\begin{rem}
Place ourselves    in the set-theoretic framework of usual cohomology of abstract groups. Then every group extension  is  admissible- provided one allows the axiom of choice, ensuring that set-theoretic splittings $\tau$ as in item (3) exist.
\end{rem}

  \begin{defi}(Direct image of an $X$-functor)\\
     Let $X \xrightarrow{f} S$ be a morphism of schemes.
Let $\Phi$ be an $X$-functor. \\Define a contravariant $S$-functor
$$\fonction{\Pi_f(\Phi)}{\{ \Sch/ S\}}{\{ \mathbf{Sets}\}}{ (T\to S) }{\Phi(X \times_S T \to X).}$$
 
 \end{defi}

 \begin{rem}\label{defifadm}
Let $X \xrightarrow{f} S$ be a morphism of schemes.
Let  \[ (\mathcal E): 1 \to A\to  B \to  C \to 1 \] be an admissible extension of  $X$-groups.
  It is straightforward to check, that \[\Pi_f(\mathcal E):  1 \to \Pi_f(A)  \to \Pi_f(B)\to \Pi_f(C) \to 1 \] is an  admissible extension of  $S$-groups.

 \end{rem}

\subsection{Affine spaces of modules.\\}
\noindent Recall some standard  and convenient notation.
\begin{defi}\label{defi PIf}
Let $ M$ be a quasi-coherent sheaf on $S$.
Define the affine space of $M$, as the functor
$$\fonctionnoname{\{ \Sch/ S\}}{ \{ \mathbf{Ab}\}}{ (Y\to S) }{H^0(Y,M \otimes_{ \mathcal O_S} \mathcal O_Y).}$$
Denote it by $\A_S(M)$, or simply $\A(M)$ if $S$ is understood. If $M$ is a vector bundle, then $\A_S(M)=\Spec(\Sym(M^\vee))$ is a smooth affine scheme over $S$.\\

\noindent Let $N_2$ be a quasi-coherent $\W_2$-module over $S$. Define its affine space as the functor 
$$\fonctionnoname{\{ \Sch/ S\}}{\{ \mathbf{Ab}\}}{(Y\to S)}{ H^0(\W_2(Y),N_2 \otimes_{ \W_2(\mathcal O_S)} \W_2(\mathcal O_Y)).}$$
Denote it by $\A_S(N_2)$, or simply $\A(N_2)$. \\If $N_2$ is a $\W_2$-bundle, then $\A_S(N_2)$ is a smooth affine scheme over $S$.

 \end{defi}

\begin{rem}
  Let $\mathbf R$ be  ring scheme over $\Z$, that is isomorphic, as a $\Z$-scheme, to an affine space $\A^n$. Then,   the second part of the preceding Definition generalises, as stated, to $\mathbf R$-modules in place of $\W_2$-modules.
\end{rem}
  
\begin{rem}
   For a quasi-coherent sheaf $M$ on $X$, the push-forward
 \[\Pi_f(\A_X(M)): \{ \Sch/ S\} \to \{ \mathbf{Ab}\}\] is given by the formula \[ (T\to S) \mapsto H^0(X \times_S T,M \otimes_{ \mathcal O_S} \mathcal O_T).\]  
\end{rem}
 
\begin{rem}
 Let $\Gamma$ be an $S$-group, and let $X \to S$ be a $\Gamma$-scheme. Consider an extension of $\Gamma$-linearised quasi-coherent  $\mathcal O_X$-modules,  of the shape \[\mathcal E: 0 \to M \to E \xrightarrow{\pi} \mathcal O_X \to 0.\] If it is geometrically trivial (=if it splits as an extension of $\mathcal O_X$-modules), it gives rise to the admissible extension of $(\Gamma,X)$-groups \[0 \to \A(M) \to \A(E) \xrightarrow{\pi} \A^1 \to 0.\]
There exists a $\Gamma$-equivariant (scheme-theoretic) section of $\pi$, iff $\mathcal E$ is trivial.

\end{rem}

   \begin{ex}
      Let  $L$ be a line bundle over $S$. The reduction sequence \[ 0 \to \A_S(L^{(1)}) \to \A_S(\W_2(L)) \to \A_S(L) \to 0\] is an admissible extension of $S$-group schemes: a  scheme-theoretic splitting is provided  by the Teichm\"uller section $\tau$. Observe that the naturality of $\tau$, ensures that it is equivariant w.r.t. to any given group action. \end{ex}

      \begin{ex}Assume that $S$ is affine, and consider an exact sequence of affine $S$-group schemes (say, w.r.t. the fppf topology) of the shape  \[ 1\to \A_S(V) \to \ast \to \ast \to 1, \] where $V$ is a vector bundle over $S$.  By Grothendieck's Theorem 90 (additive version), and vanishing of coherent cohomology over $S$, such an extension is admissible. 
  \end{ex}

\subsection{Hochschild cohomology.\\}
 
One works here  in the following context.

\begin{enumerate}
    \item{$S$ is an affine scheme, }  \item{$\Gamma \to S$ is a flat affine group scheme,}  \item{$N$ is a $(\Gamma,S)$-group functor (not necessarily an $S$-group scheme).}
\end{enumerate}

 \noindent If $N$ is commutative, the Hochschild cohomology   groups $H^*(\Gamma,N)$ are  well-defined. They  associate long exact sequences of cohomological type, to admissible extensions of commutative $\Gamma$-groups. 
They can be thought of as a sheaf-theoretic enhancement of the cohomology of finite groups. They  are defined  as  `cocycles modulo coboundaries'. This is done in the  classical  fashion: using the \textit{Hochschild complex}  \[(C^n(\Gamma,N):=N(\Gamma^n))_{n \geq 0},\] where $\Gamma^n$ is the product of $n$ copies of $\Gamma$, fibered over $S$. See \cite{J}, 4.14 for details. Observe that  we \textit{do not} use  derived functors here.\\
 As sets, $H^0(\Gamma,N)$ and $H^1(\Gamma,N)$ are still defined, without assuming $N$ commutative  (see \cite{D}, section 2), as  recalled next for the reader's convenience.

 \begin{defi}
  Let $S$ be an affine scheme, $\Gamma \to S$ be a flat affine group scheme and $N$ be a $(\Gamma,S)$-group. \\ Define $$ H^0(\Gamma, N)\subset N(S),$$ as  the subgroup consisting of $\Gamma$-invariant elements.\\ A Hochschild $1$-cocycle is a morphism of $S$-functors
$$\fonction{C}{\Gamma }{N}{\sigma}{C_\sigma,}$$
 or equivalently an element of $N(\Gamma)$, such that the cocycle relation \[ C_{\sigma \tau}=C_\sigma  .{}^\sigma C_\tau\] identically holds. Write  $Z^1(\Gamma, N)$ for the set formed by $1$-cocycles. \\Define an equivalence relation $\sim$ on $ Z^1(\Gamma, N)$,  saying that $1$-cocycles $C$ and $C'$ are equivalent (=cohomologous), if there exists $n\in N(S),$ such that, identically, \[C'_\sigma=n^{-1}. C_\sigma. {}^\sigma n  . \]  Write $$ H^1(\Gamma ,  N):= Z^1(\Gamma,  N)/\sim $$ for the factor set.\\
If $\Gamma$ acts trivially on $N$, then   $$ H^1(\Gamma ,  N)= \Hom_{S-gp}(\Gamma,  N)/{\mathrm{conjugation \; by \;}N(S)}.$$
   If $N$ is commutative, then $ Z^1(\Gamma,  N)$  has a natural abelian group structure. Cocycles cohomologous to $0$ are called coboundaries. They form a subgroup  $$ B^1(\Gamma,  N) \subset  Z^1(\Gamma,  N),$$   and $$ H^1(\Gamma ,  N)= Z^1(\Gamma,  N)/B^1(\Gamma,  N).$$
 \end{defi}

 \begin{rem}
     Equivalently (see \cite{D}, Proposition 2.2.2), a $1$-cocycle is a  (homomorphic) section, of the semi-direct product extension $$0 \to N  \to N \rtimes \Gamma   \to \Gamma \to 1.$$  Formula for the section $s$ associated to $C$: \[ s(\sigma)=(C_\sigma, \sigma).\] Cohomologous cocycles correspond to conjugate sections.
 \end{rem}

\begin{rem}
   Assume that $\Gamma$ is a finite group (seen as a constant $S$-group scheme), and  that $N$ is (the affine space of) a quasi-coherent $\mathcal O_S$-module. Then  $H^i (\Gamma,N)$ coincides with $H^i(\Gamma,H^0(S,M))$, taken in the sense of cohomology of finite groups. Here, the common notation $H^*(.,.)$  generates no  confusion.
\end{rem}

To conclude this section, recall two  classical and fundamental results.

\begin{prop}(see \cite{D}, Proposition 3.2.8)\label{HochGT}\\
 Let  $X \xrightarrow{f} S$ be a $\Gamma$-scheme, and let $M$ be a $\Gamma$-linearised quasi-coherent  $\mathcal O_X$-module. Then, the Hochschild cohomology group
\[ H^1(\Gamma,  \Pi_f(\A_X( M)))\] classifies   geometrically trivial extensions of $\Gamma$-linearised quasi-coherent  $\mathcal O_X$-modules, \[\mathcal E: 0 \to M \to E \xrightarrow{\pi} \mathcal O_X \to 0.\]

\end{prop}

For the reader's convenience, here are some (fairly classical) details.  \\Let $\mathcal E$ be as in the Proposition. Pick $e \in H^0(X,E),$ such that $\pi(e)=1$. For every $(T \to S)$, and every $\sigma\in \Gamma(T)$, one has $\pi(\sigma. e-e)=1-1=0,$ so that $$(\sigma.e-e) \in H^0(X \times_S T,M \otimes_{\mathcal O_S}\mathcal O_T).$$ The functorial association $$\fonction{C}{\Gamma }{ \Pi_f(\A( M))}{\sigma}{\sigma.e-e,}$$
  thus defines a $1$-cocycle. As one readily checks, its class in 
$ H^1(\Gamma,  \Pi_f(\A( M)))$ does not depend on the choice of $e$. In the other direction, given a $1$-cocycle in $Z^1(\Gamma,  \Pi_f(\A( M)))$ one can apply the twisting operation to the trivial extension \[ 0 \to M \to M \bigoplus \mathcal O_X \to \mathcal O_X \to 0,\] to get a geometrically split $\mathcal E$ as above. Changing the cocycle by a coboundary, yields an isomorphic extension.

\begin{rem}\label{remclassiso}
    More generally, let  $X \xrightarrow{f} S$ be a $\Gamma$-scheme, and let $\mathcal S_0$ be  `a $\Gamma$-equivariant structure over $X$' (e.g. a $\Gamma \W_2$-bundle). Then, the \textit{set} \[ H^1(\Gamma,  \Pi_f(\AAut( \mathcal S_0))) \] classifies isomorphism classes of `$\Gamma$-equivariant structures $\mathcal S$ over $X$, that are geometrically isomorphic to $\mathcal S_0$'. The proof is  the same, replacing  $(\sigma.e-e)$  by $e^{-1}  (\sigma.e)$, where $e: \mathcal S_0 \stackrel \sim \to \mathcal S$ is a geometric isomorphism. Details are left to the reader.

    \end{rem}

  \begin{rem}
    In higher cohomological degree, an important (but less definite) related notion is the Hochschild class map-- see section \ref{SeccompextH}.
  \end{rem}

\begin{prop}(`Factor systems', see \cite{D}, Proposition 2.3.6.)\label{PropFactSyst}\\
Let $N$ be a commutative  $(\Gamma,S)$-group. The Hochschild cohomology group
$ H^2(\Gamma,  N)$ classifies admissible  extensions of $S$-groups $$\mathcal E: 1 \to N \to  E \to \Gamma \to 1,$$ inducing  the prescribed action of $\Gamma$ on $N$.

\end{prop}

\section{Divided powers and Witt vectors.}\label{seceqAG}

\subsection{Divided powers.}\label{secDPA}
For a nice introduction to  divided power algebra, see \cite{Fe}. It is assumed here, that the reader is familiar with this notion.
\begin{defi}
    Let $R$ be a ring. Let $M$ be an $R$-module. For  $n \geq 1$, denote by   \[\mu (=\mu_n):  \Sym^n(M) \xrightarrow{\mu} \Gamma^n(M), \] the multiplication in the divided power algebra of $M$, given by \[x_1 \otimes \ldots \otimes x_n \mapsto [x_1]_1 \ldots  [x_n]_1.\]  Observe that $\mu(x^n)=[x]_1^n=n! [x]_n$.
\end{defi}

\noindent The following algebraic warm-up is useful.
\begin{lem}\label{LemGenSym}
    Let $n \geq 1$ be an integer, and let $R$ be a ring where $(n-1)!$ is invertible.  Let $M$ be an $R$-module.  Then following holds. \begin{enumerate}
         \item{For $i \leq n-1$, the arrow $\mu_i:\Sym^i(M) \to \Gamma ^i (M)$ is an isomorphism.}  \item{The sub-$R$-module \[\Im( \mu ) \subset \Gamma^n(M)\] is generated by elements of the shape $[x+y]_n-[x]_n-[y]_n$, for $x,y \in M$.} \item{The $R$-module $\Gamma^n(M)$ is generated by pure symbols $[x]_n$.}
         
    \end{enumerate}
   
\end{lem}

\begin{dem}
To deal with item (1), observe that the symmetrizing operator  
$$\fonctionnoname{\Gamma ^i (M)}{\Sym^i(M)}{\left[x\right]_i}{ \frac 1 {i!} x^i}$$
is well-defined for $i \leq n-1$, and provides the inverse of $\mu_i$.
Denote by $N \subset \Gamma^n(M)$ the sub-$R$-module spanned by all elements $[x+y]_n-[x]_n-[y]_n$.
   For item (2), the expansion rule for divided powers gives, for every $t \in R$, \[[ t x+y]_n-[tx]_n-[y]_n=\sum_{i=1}^{n-1} t^i[x]_i [y]_{n-i}.\]  Since $\mu_i$  is an isomorphism for $i \leq n-1$, this proves that $N \subset \Im(\mu)$. Choosing  $t=1, \ldots, n-1$ in the equality above, one gets a Vandermonde system of $n-1$ linear equations, with unknowns  the quantities $[x]_i [y]_{n-i}$, for $i=1,\ldots,n-1$. Since  $(n-1)! \in R^\times$, the corresponding  matrix is invertible.  As a consequence, for every $x,y \in M$, one has $[x]_1 [y]_{n-1} \in N$. Since  $\mu_{n-1}$  is an isomorphism,  this implies  $\Im(\mu) \subset N$, completing the proof. To prove (3), recall that,  by definition, $\Gamma^n(M)$ is generated by elements of the shape \[ [x_1]_{a_1} \ldots  [x_s]_{a_s},\] for $s \geq 1$, $a_1+ \ldots +a_s=n$ and $x_1,\ldots,x_s \in M$. By (1), one  then sees that $\Gamma^n(M)$ is generated by $\Im(\mu)$ and pure symbols, and one concludes using (2).
\end{dem}
\begin{defi}\label{defigammadiv}
    Let $R$ be a ring and let $M$ be an $R$-module. For every $a,n \geq 1$, recall  the divided power operation (\cite{Fe}, 2.2.3) 
    $$\fonctionnoname{\Gamma^a(M)}{ \Gamma^{na}(M)}{x}{ \gamma^n(x).}$$
    If $n! \in R^\times$, then $\gamma^n(x)=\frac {x^n} {n!}.$ Being functorial in $M$, it defines a polynomial law, homogeneous of degree $n$. As such, it is given by an $R$-linear map
    $$\fonctionnoname{\Gamma^n(\Gamma^a(M))}{\Gamma^{an}(M),}{\left[x\right]_n}{\gamma^n(x)}$$
that we still denote by $\gamma^n$, if no ambiguity arises.\\
 
\end{defi}
Divided powers have a dual counterpart, as follows.
\begin{defi}\label{defisigmasym}

Assumptions of Definition \ref{defigammadiv} are kept. There is a natural linear map 
$$\fonction{\sigma^n}{\Sym^{an}(M)}{\Sym^n(\Sym^a(M)),}{x_1\otimes x_2  \otimes \ldots \otimes x_{an}}{\displaystyle \sum_{\{1,\ldots,an\}=  \coprod_{j=1}^n X_j} \bigotimes_{j=1}^n( \bigotimes_{x \in X_i} x )}$$
where the sum ranges through the $\frac 1 {n!} {{na} \choose {a,a, \ldots, a}}$ unordered partitions of the set $\{1,\ldots,an\}$, into $n$ subsets of cardinality $a$. [Here `unordered' means that permuting two $X_j$'s does not change the partition.]
\end{defi}

\begin{rem}
    If the $R$-module $M$ is finite locally free, then the dual of \[\gamma^n: \Gamma^n(\Gamma^a(M))\to\Gamma^{an}(M)\]  is none other than\[\sigma^n: \Sym^{an}(M^\vee)\to  \Sym^n(\Sym^a(M^\vee)).\] The proof is standard, and it is left as an exercise.
\end{rem}
\begin{lem}\label{lemvpgamma}
    Let $R$ be a ring and let $M$ be an $R$-module. The following holds. \begin{enumerate}
        \item{For every $a,n \geq 1$,  the association 
        $$\fonctionnoname{M}{\Gamma^n(\Gamma^a(M)),}{x}{[[x]_a]_n}$$
        defines a polynomial law, homogeneous of degree $na$. [Here $\Gamma=\Gamma_R$.] As such, it is given by an $R$-linear map
        $$\begin{array}{ccc}
 \Gamma^{an}(M) & \overset{\nat}{\longrightarrow} &\Gamma^n(\Gamma^a(M)) \\
 \left[x\right]_{na} & \longmapsto & [[x]_a]_n. \end{array}$$
Then, the composite \[\Gamma^{an}(M) \xrightarrow{\nat} \Gamma^n(\Gamma^a(M)) \xrightarrow{\gamma} \Gamma^{an}(M)\] is $\frac {(na)!} {a!^n n!}\Id$.} \item{ Dually, for symmetric powers, consider the $R$-linear map 
$$\begin{array}{ccc}
 \Sym^n(\Sym^a(M)) & \overset{\mu}{\longrightarrow} &\Sym^{an}(M) \\
 f_1\otimes \ldots \otimes f_n & \longmapsto & f_1 \ldots f_n \end{array}$$
given by multiplication in the symmetric algebra $\Sym(M)$. Then, the composite \[\Sym^{an}(M) \xrightarrow{\sigma^n} \Sym^n(\Sym^a(M)) \xrightarrow{\mu} \Sym^{an}(M)\] is $\frac {(na)!} {a!^n n!}\Id$.} \item{  If  $a$ is a $p$-th power (for a prime $p$), then the integer $\frac {(na)!} {a!^n n!}$ is prime-to-$p$. If moreover $n$ is also a $p$-th power, then this integer is congruent to $1$ mod $p$. } 
    \end{enumerate}
 
\end{lem}

\begin{dem} Let us prove (1). If $R=\Q$,  then for every $b \geq 1$, the $\Q$-vector space
$\Gamma^b(M)\simeq \Sym^b(M)$ is generated by symbols $[x]_b=\frac {x^b} {b!} $, so that the composite in question is given by the formula
    \[[x]_{na}=\frac {x^{na}} {(na)!} \mapsto \frac {[x]_a^n} {n!}=\frac {x^{na}} {a!^n n!}= \frac {(na)!} {a!^n n!}[x]_{na}, \] and the conclusion readily follows from that computation. The general case is reduced to this one, via the following classical argument. Observe first that there are natural surjections (of abelian groups) \[\Gamma_{\Z}(M) \twoheadrightarrow \Gamma_R(M),\] thanks to which one reduces to the case $R=\Z$. Choose a $\Z$-linear surjection $F \to M,$ with $F$ a free $\Z$-module. This induces  surjections  $\Gamma_{\Z}(F) \to \Gamma_\Z(M),$ thanks to which we  can further assume that $M$ is free.  One can  then conclude using injectivity of the  base-change  $\Gamma(M) \to \Gamma(M) \otimes_\Z \Q= \Gamma_\Q(M\otimes_\Z \Q)$. Checking the dual formulation (2), is  straightforward using Definition \ref{defisigmasym}. Indeed, observe that \[\frac {(na)!} {a!^n n!}= \frac 1 {n!} {{na} \choose {a,a, \ldots, a}}, \] and that, in the sum expressing $\sigma^n(x_1 \otimes \ldots \otimes x_{an})$, all terms have the same image by $\mu$: namely, $x_1 \otimes \ldots \otimes x_{an}$.  For (3), by a classical theorem of Kummer, one knows that the $p$-adic valuation of a multinomial coefficient ${{a_1+\ldots a_r} \choose {a_1,a_2, \ldots, a_r}}$ is equal to the number of carry-overs,  when computing $a_1+\ldots+a_r$  in base $p$. It follows that \[ v_p({{a_1+\ldots + a_r} \choose {a_1,a_2, \ldots, a_r}})=v_p({{pa_1+\ldots  +pa_r} \choose {pa_1,pa_2, \ldots, pa_r}}),\]  for all integers $a_1,\ldots, a_r \geq 1$. This property reduces the first claim to $a=1$, where it is obvious. For the second, write $a=p^r$, $n=p^s$. Then $\frac {(p^{r+s})!} {(p^r!)^{p^s} p^s!}$ is the number of (unordered) partitions of the cyclic $p$-group $G:=\Z/p^{r+s}$, into $p^s$ subsets of cardinality $p^r$. Denote by $E$ the set formed by these partitions. Let $G$ act on $E$ by translations. It is readily checked, that the only fixed point of this action is provided by the partition formed by classes modulo the unique (cyclic) subgroup of $G$ of cardinality $p^s$. The class equation yields the result.
\end{dem}

\subsection{Witt vectors and  Teichm\"uller lifts.}\hfill\\
For an integer $r \geq 1$, we denote by $\W_r$ the $p$-typical Witt vectors of length $r$, seen as  scheme of commutative rings, defined over $\Z$.  For details, and much more on the constructions recalled below, see \cite{Bo1} and \cite{Bo2}.

\noindent To begin with, it is sufficient to know that $\W_r$ is an endofunctor of the category of  commutative  rings, such that \[\W_r(\F_p)=\Z /p^r \]and \[\W_r(\Z[\frac 1 p])=\Z[\frac 1 p]^r.\]

\noindent By gluing, $\W_r(.)$ extends to an endofunctor of the category of schemes. If $X$ is a scheme, there is a closed immersion $X \to \W_r(X)$, which is nilpotent if $X$ is a $(\Z/p^n)$-scheme, for some $n\geq 1$. On the opposite side, if $X$ is a $\Z[\frac 1 p ]$-scheme, then $\W_r(X)$ is the disjoint union of $r$ copies of $X$.\\ Thus, one can view the closed embedding $$X \hookrightarrow \W_2(X),$$ for all schemes, as a deformation of  the  usual  nilpotent  embedding $X \hookrightarrow \W_2(X) $ for $\F_p$-schemes, to the clopen embedding (on the first factor of the disjoint union)  $$X \hookrightarrow X \amalg X.$$

\subsection{Frobenius.}

\noindent Let $X$ be an $\F_p$-scheme. As would any endomorphism of $X$,  the arrow $$\frob:X \to X,$$  given by $x \mapsto x^p$ on functions, lifts to an endomorphism of $\W_r(X)$. \\ In this work, Frobenius  is  needed for  schemes  where $p$ is not a zero-divisor-- typically, projective spaces over $\Spec(\Z)$. A good reference is \cite{DK}.
\begin{mot}
    Upon replacing its source $R$ by $\W_2(R)$, the usual $\frob$ for $\F_p$-algebras, extends to all rings.  
\end{mot}
\noindent For any ring $R$, there are two natural ring homorphisms $\W_{r+1}(R) \to \W_r(R)$. The first is the reduction arrow $$\W_{r+1}(R) \xrightarrow{\rho} \W_r(R),$$ which in the coordinates $\W_r(R)=R^r$, is simply the projection $$(x_1, \ldots, x_{r+1}) \mapsto (x_1, \ldots, x_r).$$  More generally, for $1 \leq s \leq r$, the projection $$(x_1, \ldots, x_r) \mapsto (x_1, \ldots, x_s)$$ defines a ring homomorphism $$\W_r(R) \xrightarrow{\rho}  \W_s(R).$$The second is the Frobenius  $$\W_{r+1}(R) \xrightarrow{\Frob} \W_r(R),$$ given by the Witt polynomials (see \cite{DK} for details). \\Consider the natural extensions, in which surjections are ring homomorphisms, $$ 0 \to R \to\W_{r+1}(R)  \xrightarrow{\rho} \W_r(R) \to 0$$  and  $$ 0 \to \W_r(R) \xrightarrow{\Ver} \W_{r+1}(R)  \xrightarrow{\rho} R \to 0.$$ In the latter extension, provided $\W_r(R)$ is regarded as a $ \W_{r+1}(R)$-module via $\Frob$, the embedding $\Ver$ becomes $ \W_{r+1}(R)$-linear. Accordingly, it is more accurate to denote this extension  as $$ 0 \to \Frob_*( \W_r(R)) \xrightarrow{\Ver} \W_{r+1}(R)  \xrightarrow{\rho} R \to 0.$$ Observe that $ \W_{r+1}(R)$-linearity of $\Ver$,  is equivalent to the formula $$\Ver(\Frob(a)x)=a\Ver(x),$$  for all $a \in\W_{r+1}(R)  ,x \in \W_r(R)$. There is the formula $$ (\Frob \circ \Ver) = p \Id.$$ 

\begin{rem}
  Everywhere in this work,  upper-case  `$\Frob$' denotes  a Witt vector Frobenius  (over any base), whereas   lower-case `$\frob$' stands for the usual mod $p$ frobenius, given by $x \mapsto x^p$ on functions (that solely exists in characteristic $p$).
\end{rem}
\begin{rem}
Formula $(\Ver \circ \Frob = p \Id)$ is false, unless $R$ is an $\F_p$-algebra. \\ In general, by the above, $\Ver \circ \Frob $ is multiplication by $\Ver(1)$. \\ For $\F_p$-algebras, $\Frob$ equals the composite  $$\W_{r+1}(R) \xrightarrow{\rho} \W_r(R) \xrightarrow{\frob} \W_r(R).$$ 
\end{rem}

\begin{rem}
    On the opposite side, when $R$ is a $\Z[\frac 1 p]$-algebra, the product arrow $$ \W_{r+1}(R) \xrightarrow{(\rho, \Frob)} R \times \W_r(R)$$ is a ring isomorphism.

\end{rem}

A possible guideline for future research,  is suggested by the deliberately cryptic
\begin{mot}
Constructions involving ($p$-typical) Witt vectors of finite length, can actually be carried out using  iterations of $\W_2$ alone, in an appropriate way.
\end{mot}

\subsection{Equivariant $\W_r$-Modules.}\label{secgammaWmod}\hfill\\
\noindent Let $X$ be a scheme. Let $r \geq 1$. In  \cite{DCFL}, a $\W_r$-module (resp. a $\W_r$-bundle) over $X$, is defined as a quasi-coherent module (resp. a vector bundle) on the scheme $\W_r(X)$.
In this work, the $\Gamma \W_r$-modules  that occur are of scheme-theoretic nature, so that we adopt the following Definition.

\begin{defi}\label{defiGWMod}($\Gamma \W_r$-Module.)\\
Let $\Gamma \to S$ be an affine smooth group scheme, over an affine scheme $S$. 
    Let  $X \to  S$ be a $\Gamma$-scheme. Consider a commutative, affine and smooth  group scheme $\mathbf M \to X$, with the following additional data. \begin{enumerate}
     \item{ The group scheme $\mathbf M$ is endowed with the structure of a \textit {scheme in $\W_r$-modules,}  via a morphism of $X$-schemes \[\W_{r,X} \times_X \mathbf M \to \mathbf M, \] satisfying (point-wise) the usual axioms for modules over rings.}  \item{The structure  above is equipped with a $\Gamma$-linearisation. This means that an action of $\Gamma$ on $\mathbf M$ is given, over $S$ and by $\W_r$-semi-linear automorphisms.} 
 \end{enumerate}

\noindent One says that $\mathbf M$, together with the data of (1) and (2), is  a $\Gamma \W_r$-Module on $X$. 
If  $\mathbf M$ is locally free as a scheme in $\W_r$-modules, it is called a $\Gamma \W_r$-bundle  over $X$. \\
A morphism $\mathbf M \to \mathbf M'$ between  $\Gamma \W_r$-Modules, is a $\Gamma$-equivariant morphism of schemes in $\W_r$-modules.

\end{defi}
\begin{ex}
For $r \geq 1$,  let $M_r$ be a $\Gamma \W_r$-bundle over $X$. Its affine space $\A_X(M_r)$ is then a $\Gamma \W_r$-Module over $X$, in a natural way. Often, one abusively (but harmlessly) identifies $M_r$ and $\A_X(M_r)$.
\end{ex}

 \noindent If $G$ is a profinite group, define  $G \W_r$-Modules and $G \W_r$-bundles in the same way- bearing in mind that every $G$ action considered in this trilogy, factors through a suitable open normal subgroup $G_0 \subset G$.

\begin{rem}
In Definition above,  it is in general not true that the $\Gamma$-action on $X$ lifts to an action  on the $\W_r(S)$-scheme $\W_r(X)$- except if $\Gamma$ is finite étale over $S$.
\end{rem}

\subsection{Transfer w.r.t. ring schemes.}\hfill\\
 `Restricting scalars w.r.t. ring schemes'  is a widely investigated topic, that has long proved its efficiency in arithmetic geometry. In the $p$-adic setting, to the knowledge of the author, its first occurence is  Greenberg's   transform (see \cite{G}), over fields of characteristic $p$. It actually makes sense in a broader setting (see the $p$-jet spaces of \cite{Bu}).  In the present work, we  simply need to transfer  affine schemes, as follows.

 \begin{defi}\label{DefiTT}
  Let $S$ be an affine scheme.
 Let $\mathbf R$ be a ring scheme over $S$, that is isomorphic to an affine space $\A_S^r$, as an $S$-scheme.\\
    Consider the functor 
    $$\fonctionnoname{\{\Aff/S \}}{ \{\Aff/\mathbf R(S)\},}{X}{\mathbf R(X)}$$
    from affine $S$-schemes to affine $\RR(S)$-schemes. It has a right adjoint
    $$\fonctionnoname{ \{\Aff/\mathbf R(S)\}}{ \{\Aff/S\},}{Y}{\mathbf T(Y)}$$
   called the realization functor.\\ It is characterised by the following universal property. \\For any morphism    $X \to S$, there is a functorial point-wise bijection  \[ \Hom_{\Sch/S}(X,\mathbf T(Y)) \simeq \Hom_{\Sch/\mathbf R(S)}(\mathbf R(X),Y). \] 

    \end{defi}

\begin{rem}
We give the main ideas, for checking that $\mathbf T$ is well-defined. Filling in details is straightforward for an interested reader. Pick an isomorphism of $S$-schemes $\mathbf R  \stackrel  \phi \simeq    \A^r$. By a gluing argument, one may restrict w.l.o.g. to the case of affine $X$.
   Everything then becomes commutative algebra. Writing $S=\Spec(A)$, $Y=\Spec(B)$,  for an $\RR(A)$-algebra $B$,  we need to show that the functor 
   $$\fonction{\Phi_B}{\{ A-\mathbf{Alg} \}}{ \{ \mathbf{Sets} \}}{C}{\Hom_{\RR(A)-\mathbf{Alg} }(B, \RR(C))}$$
   is  represented  by an $A$-algebra.
Define \[F:=A[T_{b,i}, b \in B, i=1, 2, \ldots, r],\] the polynomial $A$-algebra on (infinitely many) indeterminates labelled by (the set) $B \times \{1, \ldots, r\}$. Consider the obvious forgetful inclusion \[  \Hom_{\RR(A)-\mathbf{Alg} }(B, \RR(C)) \subset  \Hom_{\mathbf{Sets} }(B, \RR(C))  \stackrel  \phi \simeq    \Hom_{\mathbf{Sets} }(B, C)^r = \Hom_{A-\mathbf{Alg} }(F,C),\]  where right equality  is  the universal property of the polynomial algebra $F$. \\
 One checks that there are canonical ideals  \[I^+, I^h, I^1,I^\times \subset  F,\] enjoying the following property. \\
For every $u \in \Hom_{A-\mathbf{Alg} }(F,C)$, corresponding  via the above to a function $v \in  \Hom_{\mathbf{Sets} }(B, \RR(C))$, the following equivalences hold.

\begin{enumerate}
    \item{$I^+ \subset \Ker(u)$ iff $v$ is additive (i.e. $v(x+y)=v(x)+v(y), \forall x,y \in B$).} \item{$I^h \subset \Ker(u)$ iff $v$ is homogeneous (i.e. $v( \lambda x)= \lambda v(x), \forall x \in B, \lambda \in \mathbf R(A)$).} 
      \item{$I^1 \subset \Ker(u)$ iff $v$ is unital (i.e. $v(1)=1$).}
    \item{$I^\times \subset \Ker(u)$ iff $v$ is multiplicative (i.e. $v( xy)=  v(x) v(y), \forall x,y \in B$).}
\end{enumerate}

One then checks, that the quotient $A$-algebra \[\TT(B):=F/<I^+, I^h, I^1,I^\times >\]  indeed represents the functor $\Phi_B$.
\end{rem}

\begin{rem}
    The Definition above, by generators and relations, is a common setup for  truncated $p$-jet spaces, and Weil restriction of scalars  w.r.t. a finite field extension $l/k$. Indeed, the former case is $\RR:=\W_r$, and the latter is $S=\Spec(l)$ and $\RR:=\A_k(l)$ (the affine space of the finite-dimensional $k$-vector space $l$, equipped with its natural ring scheme structure). 
\end{rem}

\begin{rem}
For our purpose, the case of an affine $Y$ (covered by Definition above) is sufficient. Let us briefly comment on the similarity between Weil restriction of scalars and  $p$-jet spaces, outside that case. The well-defineness of Weil restriction then  requires an additional assumption; typically   `every finite set of closed points of $Y$ is contained in a common open affine'.   In contrast, for $\RR=\W_r$,  and when  $p$ is nilpotent on $S$, Definition \ref{DefiTT} holds with $\Sch$ in place of $\Aff$- essentially because, for any $S$-scheme $X$, the Zariski topologies on $X$ and $\mathbf R(X)$  then coincide. Much more information on this topic (for algebraic spaces) is available in \cite{Bo2}.

\end{rem}
    \begin{defi}
In Definition \ref{DefiTT}, assume $\RR=\W_r$, for some $r \geq 1$. \\The realization functor $\TT$ is then denoted by $\TT_r$. \\Consider the (reduction) homomorphism of ring schemes $\W_r \xrightarrow{\rho} \W_1$, over $\Z$.\\
For an affine morphism $Y \xrightarrow{g} \W_r (S)$, define the cartesian diagram   \[ \xymatrix{  \rho(Y)\ar@{^{(}->}[r]^{\rho_Y} \ar[d]^{ g_S} & Y \ar[d]^g \\ S \ar@{^{(}->}[r]^-\rho &  \W_r(S), }\] whose horizontal arrows are closed immersions.
\noindent To an affine morphism $Y \xrightarrow{g} \W_r(S)$,   attach a natural morphism of $S$-schemes $$\TT_r(Y) \xrightarrow{b_Y} \rho(Y),$$  characterised point-wise  by the functorial formula, for every $(X \to S)$,\[\TT_r(Y) (X)= \Hom_{\Sch/\W_r(S)}(\W_r(X),Y) \to \Hom_{\Sch/S}(X,\rho(Y))=\rho(Y)(X)\] \[ \hspace{1.9cm}y \longmapsto \rho(y).\] Here $\rho(y)$ is defined as the factorisation of $y  \rho$, through the closed immersion $\rho_Y$, as depicted in the commutative (non-cartesian) diagram  
 \[ \xymatrix{  X\ar@{^{(}->}[r]^{\rho} \ar[d]^{ \rho(y)} & \W_r(X) \ar[d]^y \\ \rho(Y) \ar@{^{(}->}[r]^-{\rho_Y} &  Y, }\] 
\noindent [When $S=\Spec(\F_p)$, one also uses notation $\overline Y$ for $\rho(Y)$.]\\
 
 \end{defi}

 \begin{rem}
     When  $r=2$, $S=\Spec(\F_p)$, and when $g$ is smooth, Greenberg's structure theorem presents  $\TT_2(Y) \xrightarrow{b_Y} \overline Y$  as a torsor under the frobenius-twisted tangent bundle $T_{\overline Y /\Spec(\F_p)}^{(1)}$. 
 \end{rem}

\noindent Over $\Z$, one takes advantage of the obvious ring homomorphism $\Z \to  \W_r(\Z), $ leading to the following `integral' analogue of the frobenius-twisted tangent bundle.
\begin{defi}(Green bundle.)\\ \label{DefiTdagger} \noindent Let $X \xrightarrow{F} \Spec(\Z)$ be a smooth affine scheme.  Form the fibered product $$X_{ \W_r(\Z)}:= X \times_{ \Spec(\Z)}  \Spec( \W_r(\Z)) \to \Spec( \W_r(\Z)),$$ and   denote simply by   \[\TT_r X \to \Spec(\Z) ,\] the realization $ \TT_r( X_{ \W_r(\Z)}) \to \Spec(\Z) .$
Its functor of points is given by $$  \TT_r X(A)= X( \W_r(A)),$$
For any ring $A$. Consider the association,  natural in $A$,  $$X(\W_r(A)) \to X(A), $$given by composition with the ring homomorphism $\W_r(A) \xrightarrow{\rho} A$.\\
By Yoneda's Lemma, it determines
 a natural morphism $$\TT_r X \xrightarrow{b=b_{X,r}} X,$$  called the green bundle.
   By adjunction, $\Id_{\TT_r X }$ gives rise to a morphism 
    $$\mbox{ad}: \W_r(\TT_r X) \to X,  $$ such that the  following factorisation holds :  \[\xymatrix{\TT_r X \ar@/^2pc/[rr]^b \ar[r]^-\rho & \W_r(\TT_r X) \ar[r]^-{ad} & X} .\]

\noindent For a vector bundle $V$ over $X$, define the  $\W_r$-bundle over $\TT_r X$ \[\TT_r V:=\mbox{ad}^*(V). \]
 
\end{defi}

\begin{rem}
    In the Definition, the morphism $b_{X,r+1}$ clearly factors as \[ \TT_{r+1}(X) \xrightarrow{b_{X,r+1,r}}  \TT_{r}(X) \xrightarrow{b_{X,r}}  X,\] where $ b_{X,r+1,r}$ is induced point-wise by the reduction $\W_{r+1} \xrightarrow{\rho} \W_r$.
\end{rem}
\begin{ex}
    Let us work out a concrete description of the green bundle of $\GL_n$. At the light of the preceding Remark, we focus on the morphism (of $\Z$-schemes) $ b_{\GL_n,r+1,r}$, for every $r \geq 1$. There is a natural  exact sequence of affine smooth $\Z$-group schemes \[ 1 \to   \mathbf L_{r+1}(\GL_n)  \to \TT_{r+1}  (\GL_n) \xrightarrow{b_{r+1,r}} \TT_r (\GL_n) \to 1, \] whose kernel $  \mathbf L_{r+1}(\GL_n)  $ is defined point-wise as follows. Let $A$ be a ring.  Then  \[ \mathbf L_{r+1}(\GL_n) (A) \subset \mathbf M_n(A) \]  is the subset consisting of matrices $X \in \mathbf M_n(A),$ such that the matrix  \[ \Id+\Ver^r(X) \in \mathbf M_n(\W_{r+1}(A))\] is invertible.  Accordingly, the group law of $  \mathbf L_{r+1}(\GL_n) $ is given by the formula \[ X.Y=X+Y+p^r XY.\] Consider two significant cases. \begin{enumerate}
        \item{The ring  $A$ is $p$-torsion-free. Then   the invertibility condition above, is equivalent to  the simpler requirement  $ (\Id+p^rX) \in  \GL_n(A) $, and $ \mathbf L_{r+1}(\GL_n)(A)$ is just the kernel of the reduction \[ \GL_n(A) \to \GL_n(A/p^r).\]} \item{$p^r=0 \in A$. Then $\mathbf L_{r+1}(\GL_n)(A)=\mathbf M_n(A)$, with its additive structure.}
    \end{enumerate} 
\end{ex}

\begin{rem} (Inverting $p$).\\
Let $\Gamma$ be a smooth affine $\Z$-group scheme, acting on a smooth affine $\Z$-scheme $X$. Let $V$ be a $\Gamma$-bundle over $X$.  Then $\TT_r  \Gamma$ is  a smooth $\Z$-group scheme, acting on the smooth $\Z$-scheme $\TT_r X$, and the $\W_r$-bundle $\TT_r V$ is naturally a $(\TT_r  \Gamma) \W_r$-bundle.  Upon base-change to $\Z[\frac 1 p]$, the triple $(\TT_r X, \TT_r \Gamma, \TT_r V)$ becomes  isomorphic to the `direct product of $r$ copies' $(X^r,\Gamma^r, V ^r)$. The projections are induced by the `iterated Frobenii' \[\W_r \xrightarrow{\Frob^i} \W_{r-i} \xrightarrow{\rho} \W_1,\] for $i=0, \ldots, r-1$. A concrete instance of this essential fact, appears in the proof of the Smoothness Conjecture in \cite{DCF3}. In a nutshell, it simplifies otherwise intricate cohomological computations.

\end{rem}

\subsection{Teichm\"uller lift of  line bundles.\\}\label{secTeichLift}
\noindent Let $X$ be a scheme, and let $L$ be a line bundle over $X$. For $r \geq 2$, $L$ functorially extends (lifts) to a line bundle over $\W_r(X)$- its Teichm\"uller lift. We denote it by $\W_r(L)$. An introduction to the Teichm\"uller lift of line bundles can be found in  \cite{DCFL}. This is a fundamental elementary construction, obtained by applying the formalism of torsors, to the multiplicative section (of affine $\Z$-group schemes) $$ \W_1^\times=\G_m \xrightarrow{\tau} \W_r ^\times.$$ 
One can sum it up like this.

\begin{prop}\label{WittLiftTors}
    Consider the  extension of commutative  affine $\Z$-group schemes  $$0 \to \mathbf D_r \to \W_r^\times \xrightarrow{\rho} \W_1^\times(=\G_m) \to 0,$$ where $\rho$ is the natural reduction, and $\mathbf D_r$ is defined to be its kernel. It is split by the multiplicative section $\tau$.\\
     Let $L$ be a line bundle over a scheme $X$. Let $r \geq 1$ be an integer. Giving a  lift of $L$ to a $\W_r$-line bundle, is equivalent to giving a $\mathbf D_r$-torsor. The Teichm\"uller lift $\W_r(L)$ thus corresponds to the trivial torsor.
\end{prop}

\begin{dem}
    Straightforward: see \cite{DCFL}, Proposition 4.4 and Remark 4.5.
\end{dem}
\begin{rem}
Note that $\mathbf D_2 \times_\Z \F_p=\G_{a,\F_p}$. Actually, if $p>2$, $\mathbf D_r \times_\Z \F_p=\W_{r-1,\F_p}$ for all $r \geq 2$, via the $p$-adic logarithm- see \cite{DCFL}, Remark 3.5.
    
\end{rem}

\subsection{Teichm\"uller lift of vector bundles.\\}\label{secTeichLiftVB}
 Let $V$ be a vector bundle over a scheme $X$, of constant rank $d \geq 1$. In general, $V$ does not lift to a $\W_2$-bundle over $X$ (for an example where it does not, see  \cite{DCFL}, Section 7.) However, there is a functorial association \[V \to (\mathcal R \W_2(V), \tau_V)\] sending a vector bundle $V$ to a natural extension of  $\W_2$-Modules   (Definition \ref{defiGWMod}) \[\mathcal R \W_2(V): 0 \to \Frob_*(\Sym^p(V)) \to \W_2(V) \xrightarrow{\rho_V} V \to 0,\] together with a scheme-theoretic splitting of   $\rho_V$, \[ V \xrightarrow{\tau_V} \W_2(V).\] 
 \noindent How to proceed is (perhaps unsurprisingly) simple. Introduce the projective bundle \[\P(V) \xrightarrow{f} X.\] Over $\P(V)$, the twisting line bundle $\mathcal O(1)$ lifts, to $\W_2(\mathcal O(1))$. Consider the reduction sequence (of quasi-coherent  $\W_2$-modules on $\P(V)$) \[\mathcal R\W_2(\mathcal O(1)): 0 \to \Frob_*(\mathcal O(p)) \to \W_2(\mathcal O(1)) \to \mathcal O(1) \to 0.\]   Applying $f_*(.)$,  using  $f_*(\mathcal O(n))=\Sym^n(V)$  for $n \geq 0$, and $R^1f_*(\mathcal O(p))=0$, one gets  the exact sequence of   quasi-coherent  $\W_2$-modules on $X$,
 \[0 \to\Frob_*( \Sym^p(V)) \to f_*( \W_2(\mathcal O(1))) \to  V \to 0.\]  
 This construction is functorial: for any morphism $X' \xrightarrow{x} X$, set $V':=x^*(V)$, etc..., and form the exact sequence of   quasi-coherent  $\W_2$-modules on $X'$,
 \[f'_*(\mathcal R\W_2(\mathcal O_{\P(V')}(1))): 0 \to\Frob_*(\Sym^p(V')) \to f'_*( \W_2(\mathcal O_{\P(V')}(1))) \to  V' \to 0.\] 
 Define \[ \W_2(V)(X'):= f'_*( \W_2(\mathcal O_{\P(V')}(1))) .\]  The functorial extension above, defines an extension of $X$-group functors in $\W_2$-modules, that one denotes by \[\mathcal R \W_2(V): 0 \to \Frob_*(\Sym^p(V)) \to \W_2(V) \xrightarrow{\rho_V} V \to 0.\]
 It is admissible: indeed, the (functorial, sheaf-theoretic) Teichm\"uller section of $\mathcal R\W_2(\mathcal O(1))$, gives rise to a splitting $\tau_V$ of $\rho_V$, as an $X$-functor. The kernel and cokernel of $\mathcal R \W_2(V)$ are affine smooth $X$-groups schemes (as such, they should strictly speaking be denoted by $\A_X(\Sym^p(V))$  and $\A_X(V)$). Consequently, $\W_2(V) $ is an affine smooth $X$-group scheme as well, which completes the construction.
 
 \begin{rem}
     Assume that $\Gamma \to S$ is an affine smooth group scheme, over an affine base $S$, that $X \to S$ is a $\Gamma$-scheme, and that the vector bundle $V$ is $\Gamma$-linearised.  Clearly,  $\mathcal R \W_2(V)$ is  then an extension of $\Gamma\W_2$-Modules.
 \end{rem}

\begin{rem}
    The $\W_2$-Module $\W_2(V)$ is a $\W_2$-bundle, if and only if $d =1$.  In that case, i.e. if $V=L$ is a line bundle,  then $\mathcal R \W_2(V)$ is the reduction sequence of the $\W_2$-bundle $\W_2(L)$, and $\tau_V$ is its Teichm\"uller section.
 
\end{rem}
\noindent In characteristic $p$, $\mathcal R \W_2(V)$ has the following additional feature.
 \begin{lem}\label{lemkappaad}
 Let $X$ be an $\F_p$-scheme, and let $V$ be a vector bundle over $X$. The extension of $\W_2$-Modules \[\mathcal R \W_2(V): 0 \to \frob_*(\Sym^p(V)) \to \W_2(V) \xrightarrow{\rho_V} V \to 0\]
   gives rise to  a connecting arrow \[ \kappa: V \to \frob_*(\Sym^p(V)), \] defined as in section 3.5 of \cite{DCFL}. It is an $\mathcal O_X$-linear map. Via the adjunction $(\frob^*,\frob_*)$, $\kappa$ corresponds to the (divided) verschiebung of the vector bundle $V$, $\ver_V:V^{(1)} \to \Sym^p(V)$ (see section \ref{secfrobbund}).
 \end{lem}

 \begin{dem}

 Define the projective bundle $\P(V) \xrightarrow{f} X.$ 
Recall that $\mathcal R \W_2(V)$ is defined by applying $f_*(.)$ to  the reduction sequence on $\P(V)$, \[\mathcal R\W_2(\mathcal O(1)): 0 \to \frob_*(\mathcal O(p)) \to \W_2(\mathcal O(1)) \to \mathcal O(1) \to 0,\]  and that this sequence is (sheaf-theoretically) split by the Teichm\"uller section $\tau$. For a local section $s \in \mathcal O(1)$, there is the universal formula \[p\tau(s)=s^p \in \mathcal O(p),\] because  $p=0$ on $X$. Globalising yields the desired result.
 \end{dem}

\begin{rem}\label{RemVsplit}(A description of $\W_2(V)$ with coordinates).\\
    Assume that $V=\mathcal O_S^d$, with canonical basis $(e_1,\ldots,e_d)$.  Denote by $P(d)$ the set of all proper partitions $p=a_1+\ldots + a_d,$ where `proper' means that $0 \leq a_i <p$ for all $i$. There is a natural isomorphism \[\Sym^p(V) =\mathcal O_S^d \bigoplus \mathcal O_S^{P(d)}, \] where the $i$-th basis vector of $\mathcal O_S^d$ (resp. the basis vector corresponding to a proper partition $p=a_1+\ldots + a_d$) on the right side, corresponds to $e_i^p$ (resp. to $e_1^{a_1} \ldots e_d^{a_d}$) on the left side. Via this isomorphism,  $\mathcal R \W_2(V)$ reads as an extension \[ 0 \to \Frob_*(\mathcal O_S^d) \bigoplus \Frob_*( \mathcal O_S^{P(d)}) \to \W_2(V) \xrightarrow{\rho} \mathcal O_S^d \to 0.\] As such, it arises via Baer sum from two extensions of $\W_2$-Modules on $S$, \[ 0 \to \Frob_*(\mathcal O_S^d) \to \ast \to \mathcal O_S^d\to 0\] and \[ 0 \to \Frob_*(\mathcal O_S^{P(d)}) \to \ast \to \mathcal O_S^d\to 0.\] Then, the latter extension is trivial, and the former is the direct sum of $d$ copies of the reduction sequence \[\mathcal R\W_2(\mathcal O_S):  0 \to \Frob_*(\mathcal O_S) \to \W_2(\mathcal O_S) \to   \mathcal O_S\to 0.\] 
    This fact will not be used. The verification is left as an exercise to the reader.
\end{rem}

\begin{lem}\label{LemExtW2Split}
Let $V,W$ be two quasi-coherent modules on a scheme $S$. Consider an extension of $\W_2$-modules on $S$, of the shape \[\mathcal E: 0 \to \Frob_*(W) \to E \xrightarrow{f} \rho_*(V) \to 0.\] Then, $p\mathcal E$ has a canonical splitting. If moreover  the $p$-torsion of $W$ is trivial, then  \[ \Hom_{\W_2(\mathcal O_S)-Mod}(\rho_*(V),\Frob_*(W))=0,\] so that $p\mathcal E$ actually has a unique splitting.

\end{lem}
\begin{dem}
We may assume that $S=\Spec(A)$; gluing will then be automatic, from the canonical nature of the construction.  Define a function \[ \sigma: V \to E\] in the following way. For $v \in V$, pick $e \in E$ such that $f(e)=v$. Observe that $(0,1)=\Ver(1) \in \W_2(\Z)$ annihilates $\rho_*(V)$, so that $\Ver(1).e \in \Frob_*(W)$. Set \[\sigma(v):=p e -\Ver(1).e.\] We claim that $\sigma(v)$ does not depend on the choice of $e$. To check this, it suffices to prove that $p e -\Ver(1).e=0$ if $f(e)=0$. Indeed, if $f(e)=0$, then $e \in W$, so that $\Ver(1).e=\Frob(\Ver(1))e=pe$, by definition of $\Frob_*(W)$. Being well-defined, $\sigma$ is then automatically $\W_2$-linear. Clearly $f \circ \sigma=p \Id$, so that $\sigma$ provides the sought-for canonical splitting of $p\mathcal E$. Let us show that, if   $W$ is $p$-torsion-free,  then $\Hom_{\W_2(A)}( \rho_*(V),\Frob_*(W))$ vanishes. Indeed, the natural map \[W \to W[ \frac 1 p]=W \otimes_A A[\frac 1 p]\]is then injective, so that this claim can be checked upon inverting $p$, i.e. assuming $p \in A^\times$. Then, as a ring, $\mathbf R(A)=A \times A$, in such a way that  $\rho$ (resp. $\Frob$) corresponds to the first (resp. second) projection. The statement is then readily checked, for $\rho_*(V)$ (resp. $\Frob_*(W)$) is annihilated by the idempotent $(0,1)$ (resp. $(1,0)$). The last claim follows, for two splittings differ by a homomorphism 
$\rho_*(V) \to \Frob_*(W)$.
\end{dem}
\subsection{frobenius functoriality, over $\F_p$.}\hfill\\
Let $V=V_r$ be a $\W_r$-bundle, over an $\F_p$-scheme $X$. Let $m \geq 1$ be an integer. Denote by  $$\frob^m: X \to X$$ the $m$-th iterate of the (absolute) frobenius morphism of $X$.
Write $(\frob^m)^*(V)$ for the pull-back  of $V$, with respect to $\frob^m$. It is a $\W_r$-bundle. \\
 Write $(\frob^m)_*(V)$ for the push-forward  of $V$, with respect to $\frob^m$. It is a quasi-coherent $\W_r(\O_X)$-module. If $r=1$ and if $X$ is regular, then $$\frob: X \to X$$ is finite and flat, so that $(\frob^m)_*(V)$  is still a vector bundle.

\subsection{frobenius and verschiebung of  vector bundles over $\F_p$-schemes.}\label{secfrobbund}

Let $M$ be a quasi-coherent module, over an $\F_p$-scheme $X$.  The following  generalises \cite{DCFL}, Definition 6.18.

\begin{defi}
For all $a \geq 1$, there are morphisms (of vector bundles over $X$)
\begin{align*}
    \ver: \frob^*(\Sym^a(M)) &\to \Sym^{pa}(M)\\
   1 \otimes x &\longmapsto x^p
\end{align*}
and
\begin{align*}
    \frob: \Gamma^{pa}(M) &\to\frob^*(\Gamma^a(M))\\
    [v]_{pa } &\longmapsto 1 \otimes [v]_a 
\end{align*}
which we call here  the \emph{divided verschiebung} and the \emph{divided frobenius}. When generating no confusion, we may simply refer to them as  \emph{verschiebung} and \emph{frobenius}. 
\end{defi}

\begin{lem}\label{lemnatfrob}
  Recall notation of section \ref{secDPA}. The following is true. \begin{enumerate}
       \item {For every $a\geq 1$, the following diagrams commutes: \[\xymatrix{\Gamma^{pa}(M) \ar[r]^{\nat} \ar[d]^{\frob_M} &  \Gamma^p (\Gamma^a(M))\ar[d]^{\frob_{\Gamma^a(M)}} \\\Gamma^a(M^{(1)} ) \ar@{=}[r] & \Gamma^a(M)^{(1)} .}\] } 
       \item{ For every $r\geq 0$, the   following diagrams commutes: \[ \xymatrixcolsep{4pc}\xymatrix{M^{(r+1)}  \ar[d]^{\ver_M^r} \ar[r]^-{\ver_M^{r+1}} & \Sym^{p^{r+1}}(M) \ar[d]^{\sigma^p} \\\Sym^{p^r}(M^{(1)})  \ar[r]^-{\ver_{\Sym^{p^r}(M)}} &  \Sym^p(\Sym^{p^r}(M)). }\] }

   \end{enumerate}

\end{lem}

\begin{dem}
One assumes w.l.o.g. that $X=\Spec(R)$ is affine.
   It is straightforward that the two composite maps \[ \Gamma^{pa}(M) \to \Gamma^a(M)^{(1)}  \] agree on pure symbols $[x]_{pa}$.  If $R$ contains an infinite field, these generate  $\Gamma^{pa}(M)$ and the claim is proved. Otherwise, observe that the statement can be checked upon applying $(. \otimes_{\F_p} \overline \F_p)$, i.e. upon ring-change, via $R \to R \otimes_{\F_p} \overline \F_p$. This proves commutativity of (1). For (2), using Definition \ref{defisigmasym}, one computes, for $x \in M$,
   \[\sigma_p (\ver_M^{r+1}(x\otimes 1))= \sigma^p( x^{p^{r+1}})= \frac {(p^{r+1})!} {(p^r!)^p p!} x^{p^r} \otimes x^{p^r}  \otimes \ldots \otimes x^{p^r},\] which by the second part of item (3) of Lemma \ref{lemvpgamma} equals the desired \[x^{p^r}  \otimes \ldots \otimes x^{p^r}=\ver_{\Sym^{p^r}(M)}(x^{p^r} \otimes 1).\]

\end{dem}
\subsection{Frobenius functoriality,  general case.\\}
    Let $X$ be a scheme. Consider the morphism of schemes $$\Frob: \W_r(X) \to \W_{r+1}(X),$$ that functorially arises from the morphism (of ring schemes over $\Z$) $$\Frob: \W_{r+1} \to \W_r.$$ 
For   a $\W_{r+1}$-bundle   $M_{r+1}$ over $X$, \[\Frob^*(M_{r+1})=M_{r+1} \otimes_{\Frob} \W_r\]is  a  $\W_r$-bundle over $X$.\\Note the length shift (-1) in this general notion of Frobenius pull-back.\\ If $X$ is an $\F_p$-scheme, then $\Frob^*(M_{r+1})=\frob^*(M_r)$ depends only on  \[M_r:=M_{r+1} \otimes_{\rho} \W_r.\]
 For a line bundle $L$ over an arbitrary $X$, note the formula $$\Frob^*(\W_{r+1}(L))=\W_r(L^{\otimes p}).$$ 
 \subsection{The reduction sequence of a $\W_r$-bundle.}\hfill\\
Let $V_r$ be a $\W_r$-bundle, over a scheme $X$. For an integer  $1 \leq s <r$, recall the notation $V_s:=V_r \otimes_{\W_r} \W_s$. There is  the so-called   `reduction sequence' \[\mathcal R_{r,s}(V_r): 0 \to \Frob^s_*((\Frob^s)^*(V_r)) \to  V_r \to V_s \to 0,  \]  in which all three objects are actually $\W_r$-Modules in a natural way (looking at their associated affine spaces). 
\begin{ex}
  Assume that $V_r:=\W_r(L)$ for a line bundle $L$.\\ Then, its reduction sequence above reads  as $$ 0 \to \Frob^s_*(\W_{r-s}(L^{\otimes p^s})) \to \W_r(L) \to \W_s(L) \to 0.$$   It is an exact sequence of $\W_r$-Modules, split by its natural (scheme-theoretic) Teichm\"uller section.
\end{ex}

\begin{ex}
 Assume that $r=2$, $s=1$,  and  that $X$ is an $\F_p$-scheme. Then, the reduction sequence reads as  \[\mathcal R V_2: 0 \to \frob_*(\frob^*(V_1)) \to  V_2 \to V_1 \to 0.  \]  In general, it need not be scheme-theoretically split (i.e. it need not be  admissible).
\end{ex}

\section{Extensions and operations on them.}\label{secext}

\noindent Let $\Gamma \to S$ be an affine smooth group scheme, over an affine scheme $S$. \\ Let $A$, $B$  be $\Gamma \W_r$-Modules over a $\Gamma$-scheme  $X \to S$.

\begin{defi}\label{defiextGW}

An extension of $B$ by $A$, is an *admissible* exact sequence of $\Gamma \W_r$-Modules over $X$ (see  Definitions  \ref{defi admiss} and \ref{defiGWMod}), of the shape \[0 \to A \to E \to B \to 0.\] 
For any $n \geq 1$, an  $n$-extension  of $B$ by $A$, is a complex of $\Gamma \W_r$-Modules over $X$, of the shape \[\mathcal E: 0 \to A=E_0 \to E_1 \to \ldots \to E_{n-1} \to E_n=B \to 0,\]  subject to the following requirement. There are  $\Gamma \W_r$-Modules over $X$, \[A=N_0, N_1,\ldots, N_n=B,\] and extensions of  $\Gamma \W_r$-Modules (as previously defined)

\[\mathcal  N_i: 0 \to N_i \to E_i \to N_{i+1} \to 0,\] such that  $\mathcal E$ arises from these by concatenation, in the obvious way.\\
We denote by $\EExt^n_{(\Gamma \W_r,X)-Mod}(B,A)$, the collection of all $n$-extensions of $B$ by $A$. If $A$ and $B$ are $\Gamma \W_r$-bundles, denote by \[\EExt^n_{(\Gamma \W_r,X)-bun}(B,A) \subset \EExt^n_{(\Gamma \W_r,X)-Mod}(B,A) \] the sub-collection, consisting of those $\mathcal E$, where  all $E_i$'s and $N_i$'s are $\Gamma \W_r$-bundles.
\end{defi}

\begin{rems}\hfill
\begin{enumerate}\label{RemEXt}
\item{ Defining morphisms in the usual way,  $\EExt^n_{(\Gamma \W_r,X)-Mod}(B,A)$ is a category. It is in general not abelian, but one can  give a meaning to `Yoneda linked' extensions, see e.g. \cite{DF}. However, this is not needed for our purpose- whence the  deliberately vague term `collection'. }
\item{As mentionned earlier, `admissible'  is similar to `geometrically trivial'. Accordingly,  the definition above of $n$-extensions, corresponds to the \textit{strongly geometrically trivial} $n$-extensions, introduced in \cite{DCF1}, Definition 7.2. }
\item{Two cases frequently  dealt with,  are $r=n=1$ (when attempting to lift invariant global sections of $\Gamma$-vector bundles), and $r=1,n=2$ (when attempting to lift  extensions of these).}
   
     \item {When $r=1$, $S=X=\Spec(\F_p)$,  and $\Gamma =G^0$ is finite constant, what precedes boils down to    usual extensions of $\F_p[G^0]$-modules.}

\end{enumerate}
\end{rems}

\subsection{Operations.}
Extensions are subject to six frequently used functorial operations. For a more detailed exposition, see \cite{DCF1}, section 4. As one readily checks, these operations make sense here (point-wise), because the extensions of $\Gamma \W_r$-Modules considered are, by definition, admissible. Let   \[\mathcal E: 0 \to A \to * \to \ldots \to * \to B \to 0\]  be an $n$-extension of $\Gamma \W_r$-Modules over $X$.

\begin{itemize}
\item{Push-forward. If $f:A \to A'$ is an arrow in  $\{\Gamma \W_r-\mathrm{Mod} / X\}$, the push-forward \[f_*(\mathcal E): 0 \to A' \to * \to \ldots \to * \to B \to 0\] is defined in the usual fashion.}
\item{Pull-back. If $g: B' \to B$ is an arrow in  $\{ \Gamma \W_r-\mathrm{Mod} /X\}$, the pull-back \[g^*(\mathcal E): 0 \to A \to * \to \ldots \to * \to B' \to 0\] is defined in the usual fashion.}
\item{Base-change. Let $F:Y \to X$ be a (not necessarily flat) morphism of $(\Gamma,S)$-schemes. The base-change  (or pull-back, if generating no confusion)  \[F^*(\mathcal E): 0 \to F^*(A) \to * \to \ldots \to * \to F^*(B) \to 0\]   is in $\EExt^n_{(\Gamma \W_r,Y)-Mod}(F^*(B),F^*(A) )$. }
    \item {Baer sum. If \[\mathcal E_1: 0 \to A \to * \to \ldots \to * \to B \to 0\] and \[\mathcal E_2: 0 \to A \to * \to \ldots \to * \to B \to 0\] are two extensions, we can form their Baer sum   \[\mathcal E_1 + \mathcal E_2: 0 \to A \to * \to \ldots \to * \to B \to 0.\] }
    \item{Cup-product.  If \[\mathcal E_1: 0 \to A \to * \to \ldots \to * \to B \to 0\] is an $n_1$-extension, and \[\mathcal E_2: 0 \to B \to * \to \ldots \to * \to C \to 0\] is an $n_2$-extension, one can form their concatenation (or  cup-product)    \[(\mathcal E_1 \cup \mathcal E_2): 0 \to A \to * \to \ldots \to * \to C \to 0.\] It is an $(n_1+n_2)$-extension. Conversely,  a given $n$-extension breaks as a cup-product of smaller extensions, in many ways.  \begin{mot}
       Extensions are very convenient to handle cup-products. 
    \end{mot} }
    \item{Reduction. Assume that $\mathcal E$ is an extension of $\Gamma \W_r$-bundles. \\
    Then, for $1 \leq s <r$,   \[(\mathcal E\otimes_{\W_r} \W_s): 0 \to A \otimes_{\W_r} \W_s \to * \to \ldots \to * \to B \otimes_{\W_r} \W_s \to 0\] is an extension of $\Gamma \W_s$-bundles over $ X$.}
\end{itemize}
Note that push-forwards and pull-backs commute, in the sense that there is a natural isomorphism $$f_*(g^*(\mathcal E)) \stackrel \sim \to g^*(f_*(\mathcal E)).$$ Similarly, push-forward,  pull-back and base-change commute to Baer sum.

\subsection{The Hochschild class map.} \label{SeccompextH}

To fix ideas, let us consider the  case of extensions of  $\Gamma \W_r$-Modules. The  `class map' defined next, also makes sense in other setups, e.g. replacing $ \W_r$ by the ring scheme $\RR=\W_2^{[r]}$  of section  \ref{SecR2}, or dealing with strongly geometrically split extensions of  quasi-coherent $\Gamma$-linearised $\mathcal O_X$-modules.
\begin{defi}\label{defih}(Hochschild class map.)\\
    Let $S$ be an affine scheme, and let $\Gamma \to S$ be a smooth affine group scheme. Let $X \xrightarrow{f} S$ be a $\Gamma$-scheme. Let $M$ be   a $\Gamma \W_r$-Module over $X$. For every $n\geq 1$, define the  `Hochschild class map' \[ \EExt_{(\Gamma \W_r,X)-Mod}^n(\W_r,M) \xrightarrow{h} H^n(\Gamma,\Pi_f(M)) ,\] as follows. Let \[ (\mathcal E): 0 \to M=M_0 \to M_1  \to \ldots \to M_n \to \W_r\to 0 \] be an  $n$-extension of $\Gamma \W_r$-Modules over $X$.\\  By definition, it breaks into   admissible $1$-extensions,  for $i=1,\ldots, n$,  \[(\mathcal E_i): 0 \to N_{i-1} \to M_i \to N_i \to 0, \] with $N_0=M$ and $N_n=\W_r$. Denote by $\beta_i$ the Bockstein homomorphism induced by $\Pi_f(\mathcal E_i)$, in Hochschild cohomology. Successively applying these,  starting with $1 \in H^0(\Gamma, \Pi_f(\W_r))$, yields the cohomology class \[h(\mathcal  E):=(\beta_1\circ \ldots \circ \beta_{n-1} \circ\beta_n )(1) \in H^n(\Gamma,\Pi_f(M)) .\]
\end{defi}

\begin{rem}
   The Hochschild class map  is highly functorial. For instance,  one has $h(\mathcal E_1 +\mathcal E_2)=h(\mathcal E_1) +h(\mathcal E_2)$, compatibility to push-forwards (via morphisms $M \to M'$), to cup-products,  base-change and group-change... Also, when restricted to $ \EExt_{(\Gamma \W_r,X)-bun}^n(\W_r,M) $, it is compatible to the  reduction $( . \otimes_{\W_r} \W_s)$, for $s <r$.
\end{rem}
 \begin{rem}
    Let $\Gamma$ be a linear algebraic group, over a field $k$. Place ourselves in the  framework of  extensions of finite-dimensional representations of $\Gamma$ over $k$. In other words: in Definition \ref{defih}, take  $r=1$, $S=\Spec(k)$,  and $f=\Id$.   Denote by $ \Ext_{\Gamma-Mod}^n(k,M) $ the group of Yoneda-equivalence classes of extensions in  $\EExt_{\Gamma-Mod}^n(k,M) $. This group coincides with derived functor cohomology, and for all $n \geq 1$,  $h$ induces  an isomorphism, \[\Ext_{\Gamma-Mod}^n(k,M) \stackrel \sim \to  H^n(\Gamma,M) .\]  For a proof, see \cite{J},  chapter 4.
 \end{rem}

\subsection{Cohomological detox.}
\begin{mot}
  `Effectiveness of  proofs matters as much as results.'  
\end{mot}

\section{Monomial expressions.}\label{secmonexp}

The purpose of this section is to give a meaning to  `a monomial expression' applied to vector bundles over a scheme $S$. This is useful to deal with (good) filtrations later on. We give here  an elementary presentation, sufficient to deal  with applications in the paper \cite{DCF3}. For vectors spaces over fields, there is the classical notion of a polynomial functor (see  \cite{FS}), which also makes sense for quasi-coherent modules over a scheme.  In this work, we  content ourselves  with  simple  `purely multiplicative'  such functors:
 \begin{defi}
Let $n \geq 1$ be an integer. \\Over $S$, a standard monomial functor of degree $n$ is an endofunctor of the category of quasi-coherent $\mathcal O_S$-modules, of one of the following shapes. \begin{enumerate}
    \item{The $n$-th symmetric power $M \mapsto \Sym_{\mathcal O_S}^n(M)$,}  \item{The $n$-th divided power $M \mapsto \Gamma_{\mathcal O_S}^n(M)$,}  \item{The $n$-th exterior power $M \mapsto \Lambda_{\mathcal O_S}^n(M)$,}\item{If $n=p^r$ and if  $S$ is an $\F_p$-scheme, the $r$-th (absolute) frobenius twist $M \mapsto M^{(r)}=(\frob^r)^*(M)$.}
\end{enumerate}

\end{defi}

\begin{rem}
Clearly, $\Gamma^n,\Sym^n$ and $\Lambda^n$ are  defined over $\Z$ (and thus over any base), whereas $(\frob^r)^*(.)$ is defined over $\F_p$.
\end{rem}

\begin{defi}(Monomial expression.) \\ \label{DefiMonExp}
 An $S$-monomial expression in $m\geq 1$ variables, is an association $$ (V_1,\ldots, V_m) \mapsto \Theta(V_1,\ldots, V_m),$$ sending an $m$-tuple of vector bundles over $S$ to a vector bundle over $S$, that is inductively obtained by the following rules.
\begin{enumerate}
    \item{If $m =1$,  the standard monomial functors are monomial expressions.} \item{(Duality) If $m =1$, the association $(V \mapsto V^\vee)$ is a monomial expression.}  \item{(Tensor product) If $\Theta$ is a monomial expression in $m$ variables,  then \[(V_1,\ldots, V_{m+1}) \mapsto  \Theta(V_1,\ldots,V_m)\otimes  V_{m+1}\] is a  monomial expression in $(m+1)$ variables.}\item{(Substitution) If $\Theta$, $\Theta'$ are monomial expressions, in $m$ and $m'$ variables respectively, then the association \[\Theta'': (V_1,\ldots, V_{m+m'-1}) \mapsto    \Theta(\Theta'(V_1,\ldots,V_{m'}),V_{m'+1},\ldots,V_{m'+m-1})\] is a  monomial expression in  $m+m'-1$ variables.} \item{(Repetition) If $\Theta$ is a monomial expression in $(m+1)$ variables,  then \[(V_1,\ldots, V_m) \mapsto  \Theta(V_1,V_1,V_2,\ldots,V_m)\] is a monomial expression in $m$ variables.} \item{(Permutation) If $\Theta$ is a monomial expression in $m$ variables,  and if $\sigma \in S_m$  is a permutation on $m$ letters, then \[(V_1,\ldots, V_m) \mapsto  \Theta(V_{\sigma(1)},V_{\sigma(2)},V_{\sigma(3)}\ldots,V_{\sigma(m)})\] is a monomial expression in $m$ variables.}
\end{enumerate}

\noindent The degree of a monomial expression $\Theta$ is the unique $m$-tuple \[(a_1,\ldots,a_m) \in \Z^m,\] such that \[\Theta(V_1 \otimes L_1,\ldots,V_m\otimes L_m)= \Theta(V_1,\ldots,V_m) \otimes \bigotimes_{i=1}^m L_i^{\otimes a_i} ,\] whenever $L_1,\ldots, L_m$ are line bundles. The total degree of $\Theta$ is $d:=a_1+a_2+\ldots+a_m$.
\end{defi}

\begin{rem}
We justify, that the degree of a monomial expression $\Theta$ is well-defined. In the starting case where $\Theta$ is a standard monomial functor of degree $n$ (item (1) above), then clearly $\Theta$ is indeed of degree $n$ as  a monomial expression. Then, browsing through items (2) to (6), it is straightforward to show by induction that this notion makes sense in general. Let us treat the case  of item (4), and leave other cases to the reader. Assume that $\Theta$, resp. $\Theta'$, is of degree $(a_1,\ldots,a_m)$, resp. $(a'_1,\ldots,a'_{m'})$. Then it is straightforward, that the degree of $\Theta''$ is \[(a_1a'_1,a_1a'_2,\ldots,a_1a'_{m'},a_{m'+1}, a_{m'+2}, \ldots, a_{m'+m-1})\in \Z^{m'+m-1}.\]
\end{rem}
\begin{ex}
     For $m=1$, the association \[V \mapsto \End(V)=V \otimes V^\vee\]is  a monomial expression, of degree zero. Over $\F_p$, the association \[(V,W) \mapsto  \Sym^3( \Gamma^2(V) \otimes W^{(1)}) \otimes \Lambda^2(V^{(1)})\]is  a monomial expression in two variables, of degree $(6+2p,3p)$. 
\end{ex}

 \begin{rem}  
   In general,  monomial expressions are   neither covariant nor contravariant. They are, however, covariant w.r.t. isomorphisms.
\end{rem}

 \subsection{Complete filtrations and monomial expressions.}

 Monomial expressions behave well w.r.t. complete filtrations, in the sense of the next Lemma. 
 \begin{lem}\label{lemmonfil}
     Let  \[ (V_1,\ldots, V_m) \mapsto \Theta(V_1,\ldots, V_m)\] be a monomial expression of degree $(a_1,\ldots,a_m)$.\\ Assume that each vector bundle $V_i$ is equipped with a complete filtration $(V_{i,j})_{0 \leq j \leq d_i}$. [In other words, an increasing filtration by sub-bundles, such that graded pieces $L_{i,j}:=V_{i,j}/V_{i,j-1}$ are line bundles, for $j=1,\ldots, d_i$.]\\
     This data determines a complete filtration on the vector bundle $\Theta(V_1,\ldots, V_m)$, in a natural way. Graded pieces of that filtration are of the shape \[\bigotimes_{\substack{1\leq i\leq m\\ 1\leq j \leq d_i}}L_{i,j}^{\otimes \alpha_{i,j}} ,\] where the integers $\alpha_{i,j} \in \Z$ satisfy $\sum_{j=1}^{d_i} \alpha_{i,j}=a_i$, for $i=1,\ldots,m$.
 \end{lem}

 \begin{dem}
  Recall Definition \ref{DefiMonExp}. If $m=1$  and $\Theta$ is a standard monomial expression, there is a standard way to define a complete filtration on $\Theta(V)$, out of one on $V$; see Example \ref{exfilsym} for $\Theta=\Sym^n$. [There are  several natural ways, but we stick with that  one.] The same process works for $\Gamma^n$ and $\Lambda^n$. In characteristic $p$, the case of $(.)^{(m)}$ is trivial: just apply  $(.)^{(m)}$ to the filtration. It then remains to browse through the inductive process defining monomial expressions: items (2) to (6). How to proceed is straightforward. For instance, in item (3), use the tensor product complete filtration (of that on  $\Theta(V_1,\ldots,V_m)$, by that on $V_{m+1}$), defined as in section \ref{sectensgfil}. It is straightforward to treat other items.
 \end{dem}

\begin{rem}
    In Lemma \ref{lemmonfil}, `in a natural way' does not mean that the filtration is unique. For instance, a complete filtration of $V$ already determines several natural filtrations on $\Lambda^2(V)$. However, the important fact is the following. If the given filtrations on $V_1,\ldots,V_m$ are invariant w.r.t. a group action, then the filtration on $\Theta(V_1,\ldots, V_m)$ (constructed in the proof above) is invariant as well.
\end{rem}

\section{A computation in Hochschild cohomology of algebraic groups.}\label{SecHochComp}

\noindent In this  section, we work over a base field $k$.\\ The next Lemma is standard. Its short  instructive proof is included.

\begin{lem}\label{lemcoh}

Let $\G$ be a quasi-split reductive  $k$-group. Denote by $\B \subset \G $ a Borel $k$-subgroup.
Let $V$ be a  $\mathbf G$-module over $k$ (not necessarily of finite dimension).\\ For every $i \geq 0,$ the restriction arrow \[  H^i(\mathbf G,V) \xrightarrow{\Res_\B^\G} H^i(\B,V)\] is an isomorphism.
\end{lem}

\begin{dem}

If $v \in H^0(\B,V)$,  the orbit map $(g \mapsto g.v)$ factors to a morphism of $k$-varieties \[\G/ \B \to \A_k(V'),  \]  for some finite-dimensional sub-$G$-module $V' \subset V$. Having proper  connected source and affine target, this morphism is constant, implying $v \in H^0(\G,V)$. This proves the claim for $i=0$. The general case is by induction on $i$, via the usual dimension shifting argument (embedding $V$ into an injective $G$-module), and a little diagram chase left to the reader.

\end{dem}

\begin{defi}
    
    Let $d_1,\ldots,d_n \geq 1$ be  integers. Define  the algebraic $k$-group \[\mathbf G:= \GL_{d_1}\times \ldots \times \GL_{d_n},\]  and its Borel subgroup \[\mathbf B:= \B_{d_1}\times \ldots \times \B_{d_n}.\] Consider the quotient (a product of complete flag varieties) \[\Fl:=\mathbf G/\mathbf B= \Fl_{d_1}\times \ldots \times \Fl_{d_n}.\]
\end{defi}
\begin{lem}\label{lemcohvanish}(A vanishing result.)\\
    Let $n \geq 2$ be an integer, and let $E(V_1,\ldots,V_n)$ be a mononomial expression defined over $k$, homogeneous of non-zero degree. For $i=1,\ldots,n$, take $V_i:=k^{d_i}$ to be the standard representation of $\GL_{d_i}$. Consider  $E(V_1,\ldots,V_n)$ as a representation of the $k$-group $\mathbf G$ (defined above), in the natural way.\\
    Then, for all $i \geq 0$, the group \[ H^i(\mathbf G,E(V_1,\ldots,V_n))\] vanishes.

    \end{lem}
\begin{dem}

By  Lemma \ref{lemcoh}, it suffices to prove vanishing of $H^i(\mathbf B,E(V_1,\ldots,V_n))$.  Denote by  $(a_1,\ldots,a_n)$ the degree of $E(V_1,\ldots,V_n)$. We assume w.l.o.g. that $a_1 \neq 0$. We then  use the natural complete $\B$-invariant filtrations, on each $V_i$. By Lemma \ref{lemmonfil}, they yield a $\B$-invariant complete filtration on $E(V_1,\ldots,V_n)$, whose graded pieces are line bundles of the shape $L_1\otimes \ldots \otimes  L_n$, where each $L_i$ is a one-dimensional $\B_{d_i}$-representation. The center of $\B_{d_1}$, \[\mathbf  Z:= Z(\B_{d_1})\simeq \G_m\]  acts non-trivially on $L_1$, viz. by the character $(\lambda \mapsto \lambda^{a_1}$). By dévissage, it suffices to show (for each such graded piece) \[ H^i(\mathbf B,L_1\otimes \ldots \otimes  L_n)\stackrel ? =  0.\]  With the help of the spectral sequence associated to the (split) group extension \[ 1 \to \B_{d_1} \to \G \to \B_{d_2} \times \ldots \times \B_{d_n} \to 1,\] it suffices to prove \[ H^i(\B_{d_1},L_1)\stackrel ? =  0,\] for all $i \geq 0$. Consider the natural extension \[ 1 \to \U_{d_1} \to \B_{d_1} \to  \T_{d_1}=\G_m^{d_1} \to 1,\] where $\U_{d_1} \subset \B_{d_1}$ is the unipotent radical, consisting of strictly upper triangular matrices. Since algebraic tori are linearly reductive,  the natural (restriction) arrow \[ H^i(\B_{d_1},L_1) \to H^0(\T_{d_1}, H^i(\U_{d_1},L_1))\] is an isomorphism. Consider the inclusion $\mathbf Z \subset \T_{d_1}$. From the  assumption on $L_1$, it follows that the $\mathbf Z$-action on $H^i(\U_{d_1},L_1)$ occurs via the non-trivial character $(\lambda \mapsto \lambda^{a_1}$). The desired vanishing follows.

\end{dem}

  \section{Cyclotomic pairs and   smooth profinite groups.}\label{CycloSmooth}

The notion of a $(1,1)$-cyclotomic pair originates from \cite{DCF1}. Let us briefly recall it.
\begin{defi}($(1,1)$-cyclotomic pair; see \cite[§6]{DCF1})\label{DefiCyclo}\\
Let $G$ be a profinite group, and let $\Z/p^2(1)$ be a free $\Z/ p^2$-module of rank one, equipped with a continuous action of $G$. We say that the pair $(G,\Z/p^2(1))$ is $(1,1)$-cyclotomic if the following lifting property holds. \\
For all open  subgroups $H \subset G$, the natural map \[H^1(H,\Z /p^{2}(1)) \to H^1(H,\Z/p(1)) \] is surjective. 
\end{defi}

\noindent The following result is  \cite{DCF1}, Theorem A, in the particular case $n=e=1$.

 \begin{prop}\label{propcyclolift}
 Let  $(G,\Z/p^2(1))$ be a  $(1,1)$-cyclotomic pair.\\
Let $S$ be a perfect $(G,\F_p)$-scheme. Consider a $G$-line bundle $L_1$, and a geometrically split extension of $G$-linearized vector bundles over $S$, \[\mathcal E_1: 0 \to L_1(1) \to E_1  \to \mathcal O_S \to 0.\] Then,  $\mathcal E_1$ lifts to a geometrically split extension of $G\W_2$-bundles over $S$, \[\mathcal E_2: 0 \to \W_2(L_1)(1) \to E_2 \to \W_2(\mathcal O_S) \to 0.\]
\end{prop}

\section{WTF data.}\label{SecWTF}
\noindent A significant aspect of this  work, in order to build the Uplifting Pattern, is to replace  $\W_2$ by another relevant ring scheme $\mathbf R$. It should satisfy appropriate assumptions, so that  all constructions previously performed with $\W_2$-bundles, adapt to $\mathbf R$-bundles. This motivates a tentative
\begin{defi}(WTF data.) \label{defiWTFdata}\\
Let  $\mathbf R$ be a ring scheme over $\Z$, together with morphisms of $\Z$-ring schemes $$ \mathbf R \xrightarrow {\rho} \W_1(:=\A^1)$$  and $$ \mathbf R \xrightarrow {\Frob} \W_1,$$ both surjective for the fppf topology. 

\noindent Assume  given an exact sequence \[0 \to \Frob_*(\W_1) \xrightarrow {\Ver} \mathbf R \xrightarrow {\rho} \W_1 \to 0,\] where $\Frob_*(\W_1)$ equals $\W_1$ as a group scheme, and is treated as a  scheme in $\mathbf R$-modules,  by the formula $$ a.i:=\Frob(a) i,$$ for $a \in \mathbf R$  and $i \in \Frob_*(\W_1)$.\\
Assume that  the composite $$ \Frob_*(\W_1) \xrightarrow{ \Ver}  \mathbf R\xrightarrow{ \Frob} \Frob_*(\W_1)$$ equals $p \Id.$\\
Finally, assume given    a scheme-theoretic multiplicative section  of $\rho$, denoted by $$ \tau: \W_1 \to \mathbf R,$$ such that, for some integer $r\geq 0$, the following formula   identically holds: \[\Frob (\tau(x))=x^{p^{r+1}}.\]  

\noindent The quadruple $(\mathbf R, \rho,\Frob, \tau)$ is  a WTF data (Witt, Teichm\"uller, Frobenius).
\end{defi}

\begin{lem}
  Let $(\mathbf R, \rho,\Frob, \tau)$ be a WTF data. The following hold. \begin{enumerate}
      \item{The integer $r$ in Definition \ref{defiWTFdata} is unique.} \item{As a $\Z$-scheme, $\mathbf R$ is isomorphic to an affine space $\A^2$.} \item{Denote by $$\overline {\mathbf R}:=\mathbf R \times_{\Spec(\Z)} \Spec(\F_p)$$ the reduction of $\mathbf R$ to a  ring scheme over $\F_p$.  This diagram commutes: \[ \xymatrix@+2pc{\overline {\mathbf R} \ar[r]_-{\overline {\rho}} \ar@/^1pc/[rr]^{\overline {\Frob}} &\overline { \W_1} \ar[r]_-{x \mapsto x^{p^{r+1}}} & \overline { \W_1}.} \]}
  \end{enumerate}

\end{lem}

\begin{dem}
    The first item is clear. For item (2), it is straightforward to check, that the arrow
    $$\fonctionnoname{\W_1 \times_{\Z} \W_1}{ \mathbf R}{(x,y)}{\tau(x)+\Ver(y)}$$
   is an isomorphism of $\Z$-schemes. To prove the third one, observe that, upon reducing mod $p$, $\overline \Frob \circ \overline\Ver=0$. Hence, $\overline \Frob: \overline {\mathbf R} \to \overline {\W_1}$ factorises through $\overline \rho$. The statement then follows from the relation $\overline \Frob (\overline \tau(x))=x^{p^{r+1}}$.
\end{dem}
\quad\\

\noindent The next Lemma ensures that WTF data behave well, w.r.t. Zariski gluing.
\begin{lem}\label{Lem R2glues}
   Let $(\mathbf R, \rho,\Frob, \tau)$ be a WTF data. Let $A$ be a ring, and let $f \in A$. Observe that $\tau(f)$ is invertible in the ring $\mathbf R(A_f)$, so that there is a natural  ring homomorphism $$\mathbf R(A)_{\tau(f)} \to \mathbf R(A_f).$$ It is an isomorphism.\\
   If $f_1, f_2, \ldots, f_n \in A$ generate the unit ideal,  so do $\tau(f_1),\tau(f_2), \ldots, \tau(f_n) \in \mathbf R(A)$.
\end{lem}
\begin{dem}
   Identify $\mathbf R(A)$ to $A \times A$, via
   $$\begin{array}{ccc}
 A \times A & \overset{\sim}{\longrightarrow} &  \mathbf R(A) \\
 (x,y) & \longmapsto & \tau(x)+\Ver(y). \end{array}$$   
   Observe that, for all $x,y \in A$, $$ \tau(f) (x,y)=(fx,f^{p^{r+1}}y) \in \mathbf R(A).$$One can thus define an arrow $$\mathbf R(A_f)\longrightarrow \mathbf R(A)_{\tau(f)}$$  $$\hspace{-2.15cm}(\frac x {f^a}, \frac y {f^b})=(\frac {f^{N-a} x} {f^N}, \frac {f^{Np^{r+1}-b} y} {f^{Np^{r+1}}}) \longmapsto \frac {(f^{N-a} x,f^{Np^{r+1}-b} y)} {\tau(f)^N}.$$ [This formula does not depend on the choice of $N > \mathrm{max}(a,b)$.]\\ As one readily checks,  it  defines the inverse of the homomorphism in  the first part of the Lemma. For the second part, observe that, if $f_1, f_2, \ldots, f_n \in A$ generate the unit ideal,  so do $f_1^{p^{r+1}},f_2^{p^{r+1}}, \ldots, f_n^{p^{r+1}} \in A$. The statement is then a straightforward two-step dévissage, using the formula $ \tau(f) (x,y)=(fx,f^{p^{r+1}}y)$.  
\end{dem}

\begin{rem}
     Let $(\mathbf R, \rho,\Frob, \tau)$ be a WTF data. The morphism of $\F_p$-ring schemes \[ \overline { \W_1} \xrightarrow{x \mapsto x^{p^{r+1}}} \overline { \W_1}\] will be denoted by $\frob_\RR$, or by $\frob^{r+1}$. It is the `relevant' mod $p$ frobenius arrow, when working with the ring scheme $\RR$ in place of $\W_2$.
\end{rem}

\subsection{Teichm\"uller lift of vector bundles, w.r.t. $\RR$.\\}\label{secTeichLiftVBR}

The construction of section \ref{secTeichLiftVB} adapts \textit{verbatim} to the broader setting of a WTF data $(\mathbf R, \rho,\Frob, \tau)$.  Results are summarized below.\\
For  a vector bundle $V$ over a scheme $X$, introduce its projective bundle \[\P( V) \xrightarrow{f_V} X.\]  The  functorial association \[V \to \RR(V):=(f_V)_*(\RR(\mathcal O(1))),\] sends a vector bundle $V$ to an $\RR$-Module $\RR(V)$, together with a natural (admissible) extension of  $\RR$-Modules \[\mathcal R \RR(V): 0 \to \Frob_*(\Sym^{p^{r+1}}(V)) \to \RR(V) \xrightarrow{\rho_V} V \to 0.\] If  $V=L$ is a line bundle,   $\mathcal R \RR(V)$ is the reduction sequence of the $\RR$-bundle $\RR(L)$.

\subsection{The multiplicative group scheme $\mathbf R^\times$.}
\begin{defi}\label{defiGam}
Taking invertible elements in the ring scheme $\mathbf R$, yields the multiplicative group scheme $$\mathbf R^\times  \to \Spec(\Z).$$ It fits into an exact sequence of commutative, affine  and smooth $\Z$-group schemes
\[ 1 \to \G_{a/m} \to \mathbf R^\times  \xrightarrow{\rho} \G_m \to 1,\]
split by the multiplicative section $\tau$. The notation $\G_{a/m}$ for its kernel, is justified by the fact that  it is a deformation of the additive group to the multiplicative group. Indeed, its generic fiber is multiplicative: $$\G_{a/m}\times_{\Spec(\Z)}{\Spec(\Q)}=\G_m \times_{\Spec(\Z)}{\Spec(\Q)},$$ whereas  its special fiber is additive:  $$\G_{a/m}\times_{\Spec(\Z)}{\Spec(\F_p)}=\G_a \times_{\Spec(\Z)}{\Spec(\F_p)}.$$

\end{defi}
\begin{prop}(Lifts of line bundles, seen as $\G_{a/m}$-torsors.)\label{propgamtors}\\
 Let $L$ be a line bundle over a scheme $X$.  Giving a  lift of $L$ to an $\mathbf R$-line bundle, is equivalent to giving a $\G_{a/m}$-torsor. The Teichm\"uller lift $\mathbf R(L)$ thus corresponds to the trivial torsor.
\end{prop}

\begin{dem}
    Same as Proposition \ref{WittLiftTors}, using the ring scheme $\mathbf R$ in place of $\W_2$.
\end{dem}
\subsection{ The functor $ X \mapsto \mathbf R(X)$, for schemes.\\}

\begin{defi} \label{defiRsch}

    Let $X$ be a scheme. If $X=\Spec(A)$ is affine, define $$\mathbf R(X):=\Spec(\mathbf R(A)).$$ In the general case, cover $X$ by open affines $X_i=\Spec(A_i)$. The affine schemes $\mathbf R(X_i)$ then glue naturally. Checking this is routine, with the help of Lemma \ref{Lem R2glues}. This gives rise to the scheme $\mathbf R(X)$.\\
  One  denotes by  $$\rho_\RR: X \to \mathbf R(X),$$ resp.  $$\Frob_\RR: X \to \mathbf R(X),$$the closed immersion, resp. the morphism, induced by $\rho_\RR$, resp. by $\Frob_\RR$. When no confusion is possible, the subscript $\RR$ may be dropped from notation.\\  For a given $\mathbf R$-bundle $V_\RR$ over $X$,  use the notation $$ V:=\rho^*(V_\RR)$$ and $$ V^{[1]}:=\Frob^*(V_\RR).$$ These are vector bundles over $X$. Strictly speaking, $V^{[1]}$  should be denoted by $ V_\RR^{[1]}$, as it actually depends on $V_\RR$ and not solely on $V$ (unless $p=0$ on 
    $X$).  Meanwhile, it will always be clear in practice, which $V_\RR$ was used to define  $V^{[1]}$. \\
 Accordingly, if $X$ is an $\F_p$-scheme, set $$ \frob_\RR=\frob^{(r+1)}: X \to X.$$  Here the subscript  $\RR$ is  kept, for an obvious reason.\\ For a vector bundle $V$ over an $\F_p$-scheme $X$,  set $$ V^{[1]}:=(\frob^{r+1})^*(V),$$ which is consistent with the above.

\end{defi}

\section{The ring scheme $\W_2^{[r]}$  over $\Z$, and associated WTF data.}\label{SecR2}

\noindent This ring scheme  is  used in the construction of  the Uplifting Scheme. Besides truncated Witt vectors themselves, it is  the only  ring scheme featured in this text. It depends on a given integer $r \geq 0$. When $r=0$,  $\W_2^{[r]}=\W_2$.
\begin{mot}
    Over an $\F_p$-base, lifting a given vector bundle $V$ to an  $\W_2^{[r]}$-bundle, is equivalent to lifting its $r$-th (absolute) frobenius twist $V^{(r)}$, to a $\W_2$-bundle.

\end{mot}

\noindent Over $\F_p$, $ \W_2^{[r]}$ can be defined by pushing-forward the reduction sequence of $\W_2$, via the homomorphism $\frob^r$ (see Lemma \ref{LemR2overFp} for a precise formulation). However, in this work, we  crucially need  $ \W_2^{[r]}$ to fit into a WTF data; in particular, to be defined over $\Z$. Its construction  then exploits truncated Witt vectors of size $(r+2)$. 
\begin{defi}\label{DefiW2r} 
    Let $r \geq 0$ be an integer. Consider the exact sequence of $\W_{r+2}$-Modules over $\Z$,\[ \mathcal E \W_{r+2}: 0 \to \Frob_*(\W_{r+1})  \to  \W_{r+2} \to \W_1 \to 0,\] and the homomorphism of ring schemes $$\Frob^r: \W_{r+1} \to \W_1,$$   given by the Witt polynomial $$(X_0,X_1, \ldots, X_r) \mapsto X_0^{p^r} +pX_1^{p^{r-1}}+ \ldots+ p^r X_r. $$  Define a ring scheme $$\W_2^{[r]}  \to \Spec(\Z) $$ by forming the push-forward diagram  $$ \xymatrix{ \mathcal E \W_{r+2}:0 \ar[r] & \W_{r+1}  \ar[r]^-{\Ver} \ar[d]^{\Frob^r} & \W_{r+2} \ar[r]^\rho \ar[d]^\pi &  \W_1 \ar[r] \ar@{=}[d]  & 0 \\ \mathcal E \W_2^{[r]}: 0 \ar[r] & \W_1  \ar[r] &  \W_2^{[r]} \ar[r]^{\rho^{[r]}} & \W_1 \ar[r]   & 0. } $$

   \noindent On the bottom line, the arrow $\W_1  \hookrightarrow \W_2^{[r]}$ is also denoted by $\Ver$.\\
The Teichm\"uller section $\tau$ of $\rho$ gives rise to    a section of $\rho^{[r]}$, $$ \tau^{[r]}: \W_1 \to \W_2^{[r]},$$  also scheme-theoretic and multiplicative, and referred to as the  Teichm\"uller section.  \\
   Recall that the composite $$\W_{r+1}\xrightarrow{ \Ver} \W_{r+2} \xrightarrow{ \Frob }\W_{r+1}  $$ is $p \Id$. It follows that $$\Frob^{r+1}: \W_{r+2} \to \W_1 $$  passes to the quotient via $\pi$, to  a natural homomorphism of ring schemes $$\Frob^{[r+1]}: \W_2^{[r]} \to \W_1, $$  and that  the composite $$\W_1 \xrightarrow{\Ver}  \W_2^{[r]}\xrightarrow{\Frob^{[r+1]} }\W_1  $$ is $p \Id$. Thus, the ideal structure on $\W_1 \subset \W_2^{[r]} $ is given by the formula $$ x.i= \Frob^{[r+1]}(x)i,$$ where multiplication on the right side, is that of the ring scheme $\W_1$. Taking into account this ideal structure,  it is accurate to   denote $ \mathcal E \W_2^{[r]}$ by $$ \mathcal E \W_2^{[r]}: 0 \to \Frob^{[r+1]}_*(\W_1)  \to  \W_2^{[r]} \to \W_1 \to 0.$$

\end{defi}

\noindent A word of explanation is perhaps needed, to justify representability of $ \mathcal E \W_2^{[r]}$.   To start with,  observe that $\Frob^r: \W_{r+1} \to \W_1 $  is an fppf morphism of $\Z$-schemes- hence surjective for the fppf topology. Consider the ideal $$\Ker(\Frob^r) \subset \W_{r+1},$$ as an fppf sheaf (actually represented by an an affine $\Z$-group scheme).  Define $\W_2^{[r]} $  to be the fppf quotient of $\W_{r+2}$, by its fppf sheaf of ideals $\Ver(\Ker(\Frob^r))$. This way, one defines $ \mathcal E \W_2^{[r]}$, as an exact sequence   of fppf  sheaves (of abelian groups) on $\Sch/ \Z$. Considered as a surjection of fppf sheaves,  $\rho^{[r]}$ has a section, given by $\tau^{[r]}$. It is then  straightforward to check  that  $\W_2^{[r]}$ is indeed representable over $\Z$, by the affine space $\A^2=\W_1 \times \W_1$.

\begin{defi}(The ring scheme $\RR= \W_2^{[r]}$ and related content.) \label{defi R2sch}\\
    For convenience, from now on, denote \[ \W_2^{[r]}, \rho^{[r]}, \Frob^{[r+1]}, \tau^{[r+1]} \] simply  by \[ \RR,    \rho_\RR, \Frob_\RR, \tau_\RR.\]The dependence on $r$ is implicit throughout.  If there is no risk of confusion with  the `usual' Witt vector Frobenius, we may drop the subscript $\RR$ from notation.\\Besides Witt vectors $\W_n$ themselves, it is the only ring scheme  used in this work.\\  The quadruple $(\RR, \rho_\RR,\Frob_\RR, \tau_\RR)$ is a WTF data (see Definition \ref{defiWTFdata}.)\\
    Let $X$ be a scheme. 
  One  denotes by  $$\rho: X \to \RR(X),$$ resp.  $$\Frob: X \to \RR(X),$$the closed immersion, resp. the morphism, induced by $\rho$, resp. by $\Frob$.\\  For a given $\RR$-bundle $V_\RR$ over $X$,  use the notation $$ V:=\rho^*(V_\RR)$$ and $$ V^{[1]}:=\Frob_\RR^*(V_\RR).$$ These are vector bundles over $X$. Observe that $ V^{[1]}$ actually depends on $V_\RR$ and not solely on $V$- unless $p=0$ on 
    $X$.\\
 Accordingly, if $X$ is an $\F_p$-scheme, put $$ \frob_\RR:=\frob^{r+1}: X \to X.$$ It is given by raising functions to their $p^{r+1}$-th power.
\end{defi}

\begin{lem}\label{Lembrpr}
    Let $V_{r+2}$ be a $\W_{r+2}$-bundle, over a scheme $X$. Consider the natural surjection (of ring schemes over $\Z$) $\W_{r+2} \xrightarrow{\pi} \RR.$ In Definition above, take $V_\RR:=V_{r+2} \otimes_\pi \RR$. Then, one has \[V^{[1]}=(\Frob^{r+1})^*(V_{r+2})=V^{(r+1)}.\]
\end{lem}
\begin{dem}
    This is  because, by definition, the composite $ \W_{r+2} \xrightarrow{\pi} \RR \xrightarrow{\Frob_\RR} \W_1$ equals  $ \W_{r+2} \xrightarrow{\Frob^{r+1}} \W_1$ (see Definition \ref{DefiW2r}).
\end{dem}

\quad\\
\noindent Base-changed to $\F_p$, the ring scheme $\RR$ can be described  using $\W_2$ alone. 

\begin{lem}\label{LemR2overFp}
    Denote by $$\overline{\RR}  \to \Spec(\F_p) $$ the base-change of $\RR$, to a ring scheme over $\F_p$.\\ It naturally fits into the push-forward diagram $$ \xymatrix{ \mathcal R \overline{\W}_2:0 \ar[r] & \frob_*(\overline{\W}_1)  \ar[r]^-{\Ver} \ar[d]^{\frob^r} & \overline{\W}_2 \ar[r]^\rho \ar[d] &  \overline{\W}_1 \ar[r] \ar@{=}[d]  & 0 \\  \mathcal R \overline{\RR}:  0 \ar[r] & \frob^{r+1}_*(\overline{\W}_1)  \ar[r] &  \overline{\RR} \ar[r] & \overline{\W}_1 \ar[r]   & 0. } $$  Furthermore, the arrow $${ \overline \Frob}^{[r+1]}: \overline{\RR} \to  \overline{\W}_1$$ equals the composite $$  \overline{\RR} \xrightarrow{ \rho}  \overline{\W}_1 \xrightarrow{\frob^{r+1}}  \overline{\W}_1.$$
\end{lem}

\begin{dem}
    Straighforward verification, left to the reader.
\end{dem}

\begin{lem}
    The ring $\RR(\F_p)$ is  $\Z/p^2$.
\end{lem}
\begin{dem}
Evaluate the diagram in the preceding Lemma at $\F_p$-rational points- remembering that $\W_2(\F_p)=\Z/p^2$.
\end{dem}
\begin{lem}\label{LemCar}
   For each $n \geq 1$, the ring $\RR(\Z/p^n)$ is a $(\Z/p^{n+1})$-algebra, in which $p^n \neq 0$. Hence, if $X$ is a $\Z/p^n$-scheme, $\RR(X)$ is a $\Z/p^{n+1}$-scheme.
\end{lem}

\begin{dem}
  By definition of $\RR$, the statement is equivalent to the following fact.\\
  For every $r \geq 0$ and $n \geq 1 $, the element  $1 \in \W_{r+1}(\Z/p^n)$ has (additive) order $p^{r+n}$. \\ Let us solve this exercise, by induction on $n+r$. For the starting cases, observe that  $r=0$ is trivial, and that $n=1$ is standard. In the ring $\W_{r+1}(\Z/p^n)$, one has  $$ p=(p, 0, 0, 0)+ (0, x_1, \ldots, x_r)=\tau(p)+ \Ver(X),$$ for some $$X=(x_1, \ldots, x_r) \in \W_r(\Z/p^n). $$ Since $p \neq 0 \in \W_2(\F_p)$, one must have $x_1 \neq (p\Z/ p^n\Z)$, so that $X$, regarded as an element of the ring $\W_r(\Z/p^n)$, is invertible. Compute, in the ring $\W_{r+1}(\Z/p^n)$: \[ p^2=(\tau(p)+ \Ver(X))^2=\tau(p)^2+2\tau(p)\Ver(X)+ \Ver(X)^2\]
   \[=\tau(p^2)+2 \Ver(p^p X)+ \Ver( p X^2),\] using multiplicativity of $\tau$, and $(\Frob\circ  \Ver)=p \Id $. By the induction hypothesis, one knows that $1 \in \W_{r+1}(\Z/p^{n-1})$ has order $p^{r+n-1}$.  Thus,   $\tau(p^2) \in \W_{r+1}(\Z/p^{n-1})$ has order $p^{r+n-3}$. Similarly, one knows that $1 \in \W_r(\Z/p^n)$ has order $p^{r+n-1}$. Since $X^2 \in \W_r(\Z/p^n) $ is invertible,  is follows that  the order of $p X^2 \in \W_r(\Z/p^n)$  is $p^{r+n-2}$.  For the same reason (and because $p \geq 2$), $p^p X \in \W_r(\Z/p^n) $ is  of order $\leq p^{r+n-3}$. Thus, the order of $p^2 \in \W_{r+1}(\Z/p^n)$ is $p^{r+n-2}$. Equivalently,  $1 \in \W_{r+1}(\Z/p^n)$ has order $p^{r+n}$, completing the induction step.
\end{dem}

\begin{exo}
    Give a precise description of $\RR(\Z/p^2)$, and show there exists a unique ring homomorphism $\RR(\Z/p^2) \to \Z/p^3$.
\end{exo}
\begin{rem}
I don't know any example of a `meaningful' WTF data, besides $\RR=\W_2^{[r]}$ (that's one for each $r\geq 0$). The reader may want to investigate this. 
\end{rem}

 \section{Splitting schemes and good filtrations.}\label{SecAff}
In this section, $S$ is a scheme. 
 \subsection{Affine spaces.} \label{DefiAffSpace}\hfill\\
 Let  $M$ be a quasi-coherent $\mathcal O_S$-module.
Recall that the affine space  $$\mathbb A_{S}(M)\to S$$ is  the contravariant functor, from the category of $S$-schemes to that of  abelian groups, defined  by the formula \[ ( T \to S) \mapsto H^0( T,M \otimes_{\O_{S}}\O_{T}).\] [It is the functor $\underline M$ used to define \textit{polynomial laws} in \cite{Fe}, Section 2.2.1.]\\
If $M$ is a vector bundle, $\mathbb A_{S}(M)$ is  represented by the smooth affine morphism
 \[\Spec(\bigoplus_{n=0}^\infty\Sym^n (M^\vee)) \to S.\]

 \begin{rem}\label{RemAffInd}  Assume that  $M$ is equipped  with an increasing  filtration $$M =\varinjlim_{j \in \N} M_j,$$ with $M_0=0$. In that case, we have  \[\mathbb A_{S}(M)=\varinjlim_{j \in \N}\mathbb A_{S}(M_j) . \]
 If moreover all quotients  $M_{j+1}/M_j$ are vector bundles, then  so are the $M_j$'s, and $\mathbb A_{S}(M)$ is an  ind-$S$-scheme, in the usual sense: for $j < j',$ the natural arrow $$\mathbb A_{S}(M_j) \to \mathbb A_{S}(M_{j'})$$ is a  closed immersion of smooth $S$-schemes.
 \end{rem}

 \subsection{Splitting schemes.}\label{SplittingSchemes}

The following notion was introduced in \cite[§4]{DCF1}- without equivariance, that is added for free.
\begin{defi}(Splitting scheme.)\\\label{GAffRep}
 Let $\Gamma \to S$ be a (smooth affine) group scheme, and let $X \to S$ be a $\Gamma$-scheme. Consider an extension of $\Gamma$-vector bundles over $ X$, \[\mathcal E: 0 \to V \stackrel i \to E \stackrel \pi \to \O_{ S} \to 0.\] Then, the $ (\Gamma,S)$-scheme of sections of $\pi$ is denoted by $$g:\mathbb S(\mathcal E) \to  X.$$It is a $\Gamma$-affine subspace of $\mathbb A_X(E )$,  and a $\Gamma$-equivariant torsor under $\A_X(V)$.\\As such, it is the $\Spec$ of the $ \Gamma \mathcal O_X$-algebra \[\lim\limits_{\longrightarrow}(\Sym^{n}(E^\vee)), \] where  limit is taken with respect to the  injections of the exact sequences  \[0 \to   \Sym^{n}(E^\vee) \stackrel { \times \pi^\vee} \to  \Sym^{n+1}(E^\vee) \stackrel {\Sym^{n+1}(i^\vee)}\to \Sym^{n+1}(V^\vee) \to 0, \]which are the symmetric powers of the dual extension \[\mathcal E^\vee: 0 \to \O_{\overline S} \stackrel {\pi^\vee} \to E^\vee \stackrel {i^\vee} \to V^\vee \to 0.\]
This description yields a natural ($\Gamma$-equivariant) filtration on the quasi-coherent $\O_X$-module $g_*(\O_{\mathbb S(\mathcal E)})$, by the  sub-vector bundles $\Sym^{n}(E^\vee)$. It is indexed by the well-ordered set $\N$. Its $n$-th graded piece is  the vector bundle $\Sym^n(V^ \vee)$. 
\end{defi}

\subsection{good filtrations.}

Filtrations (invariant w.r.t. some group action) are a main tool in this work.  They are natural and enjoy nice properties, which facilitate dévissage arguments. They  are increasing filtrations labelled by   well-ordered sets- typically $\N^d$ with lexicographic order. 
\begin{defi}(Good filtrations, well-filtered morphisms)\\
Let  $M$ be a quasi-coherent $S$-module.\\
A good filtration on $M$ is the data of 
\begin{itemize}
    \item { A well-ordered set $J$, whose least element we denote by $0$.}
    \item{An increasing filtration $(M_j)_{j\in J}$ of $M$, by quasi-coherent sub-modules, whose graded pieces $$F_j:=M_j/\sum_{i < j} M_{i}$$ are vector bundles, for all $j \in J$. If all $F_j$'s are line bundles, the filtration is called complete.}
\end{itemize}
Let $g:S' \to S$ be an affine morphism. Let $(M_j)_{j\in J}$ be a good filtration of the quasi-coherent $S$-module $g_*(\mathcal O_{S'})$, with first step $$M_0=\mathcal O_S \subset g_*(\mathcal O_{S'}).$$ We then say that $g$ is well-filtered, w.r.t. the good filtration  $(M_j)_{j\in J}$.
\end{defi}

\begin{rem}
    Assume that $\Gamma \to S$ is an affine group scheme, that $M$ is $\Gamma$-linearised, and that we need  to show vanishing of $H^i(\Gamma,M)$. By dévissage on a $\Gamma$-invariant good filtration as above, it then suffices to prove  $H^i(\Gamma,F_j)=0$, for all $j \in J$. This is straightforward to check, using that $J$ is well-ordered. 
\end{rem}

\subsection{Tensor product of good filtrations. }\label{sectensgfil} \hfill\\
Let $M,M'$ be quasi-coherent modules over $S$, respectively equipped with good filtrations  $(M_{j})_{j\in J}$ and $(M'_{j'})_{j'\in J'}$.  Put $$M_{j, j'}:= \sum_{ (i,i')\leq (j,j') } M_i \otimes M'_{i'}.$$Then, $(M_{j, j'})_{(j,j')\in J \times J'}$ is a good filtration of $M \otimes M'$, for the lexicographic order on $J \times J'$. Its graded pieces are the vector bundles $$F_{j,j'}:=(F_j \otimes F_{j'})_{(j,j')\in J \times J'}. $$This extends to tensor products of three or more good filtrations of quasi-coherent modules. This construction depends on the choice of an order of the factors.

 \begin{rem}(Fibered product of well-filtered morphisms.)\\
  Assume that $T_1 \to S$ and $ T_2 \to S$ are two well-filtered  morphisms.  Their fibered product \[ T_1 \times_S T_2 \to S\] is then naturally   well-filtered, via  the tensor product of good filtrations above. \\This  extends to  fibered products of three or more factors.
 \end{rem}

 \subsection{Examples of good filtrations.}
 \begin{ex}(A typical good filtration, on symmetric powers.)\label{exfilsym}\\
 Let $S$ be a scheme.  Let $V$ be a $d$-dimensional vector bundle over $S$, equipped with a complete filtration by sub-bundles, \[\nabla: 0 \subset  V_1 \subset \ldots \subset V_d=V.\] Denote by $L_i:=V_i/V_{i+1}$ the line bundles forming its graded pieces. For  every  $n \geq 1$, the vector bundle $  \Sym^n(V)$ then has a natural good filtration $(M_{j})_{j\in J}$. Here  $$J:=\{( a_1,\ldots,a_{d}) \in \N^{d}, \; \sum a_i=n\}$$ is endowed with the lexicographic order, and \[M_{( a_1,\ldots,a_{d})}:=\mathrm{Span}( v_1\ldots v_d, \; v_i \in  V_{b_i},\; (b_1,\ldots,b_{d}) \leq (a_1,\ldots,a_{d})). \] Graded pieces of that filtration are $F_{(a_1,\ldots,a_{d})}=L_1^{\otimes a_1} \otimes \ldots \otimes  L_d^{\otimes a_d}$.

 \end{ex}

 \begin{ex}
Splitting schemes of extensions of vector bundles (see Definition \ref{GAffRep}) are  archetypes of  well-filtered  morphisms. 
 \end{ex}

\section{Lifting frobenius: the uplifting scheme.}\label{SecLiftFrob}
\begin{center}
  \textit{In this section, $S$ is a scheme such that $p$ is nowhere a zero-divisor on $S$.}
\end{center}
[Equivalently, for every open $U \subset S$, the ring $\mathcal O_S(U)$ is $p$-torsion-free.] \\
 In typical applications later on, $S$ is a smooth scheme over $\Z$. Set $\overline S:=S \times_\Z \F_p$.
\subsection{On equivariance.}

In the exposition that follows, constructions  are, once more, entirely canonical. To keep notation light, their good behaviour w.r.t. equivariance is thus implicit.

\subsection{The uplifting scheme, for $\W_2$ coefficients.\\}

Let $V$ be a vector bundle over $S$, of constant rank $d\geq 1$. As usual, set  $\overline V:=V \otimes_{\Z} \F_p$. Recall the (divided) verschiebung homomorphism for vector bundles over $\overline S$ (see section \ref{secfrobbund}), \[ \ver_{\overline V}: \overline V^{(1)} \to \Sym^p(\overline V).\]
Our next goal, is to produce a natural (affine and smooth) base-change \[\mathbf U \to S,\] such that, over $\mathbf U$,  $\ver_{\overline V}$ acquires an integral lift. Equivalenty, the frobenius $\frob_{\overline V^\vee}$ acquires an integral lift. For this assertion to make sense, it must be the case that $\overline V^{(1)}$ itself lifts, to a vector bundle over $S$. To ensure this, \begin{center}
    \textit{one  assumes  given a lift of $V$, to a $\W_2$-bundle $V_2$ over $S$. }
\end{center}Then, the assertion is precisely stated like this.\\ Over $\mathbf U$, $\ver_{\overline V}$ lifts, to a homomorphism of vector bundles  \[ \Ver_V:  V^{(1)} \to \Sym^p(V),\] and $\mathbf U$
 is universal w.r.t. this property.\\
To proceed, consider $\ver_{\overline V }$ as a global section $$\ver_{\overline V }\in H^0(\overline S,   \overline  V^{\vee (1)} \otimes  \Sym^p(\overline V)).$$ Since $p$ is nowhere a zero-divisor on $S$, there is an exact sequence of $\mathcal O_S$-modules \[ 0 \to V^{\vee (1)} \otimes  \Sym^p(V) \xrightarrow{p \Id} V^{\vee (1)} \otimes  \Sym^p(V)  \to   \overline  V^{\vee (1)} \otimes  \Sym^p(\overline V)  \to 0.\] From there, one sees that the existence of $\Ver_V$ is obstructed by a natural class $$\obs(\Ver_V) \in H^1(S,    V^{\vee (1)} \otimes  \Sym^p(V)).$$ 

\begin{defi}\label{defiUW2}
Define the uplifting scheme of $V_2$, as the splitting scheme (see Definition \ref{GAffRep}) \[\mathbf U(V_2):=\mathbb S(\obs(\Ver_V)) \to S.\] Denote by  \[ \Ver_V:  V^{(1)} \to \Sym^p(V)\] the universal  Verschiebung, lifting $\ver_{\overline V}$ over $ \mathbf U =\mathbf U(V_2)$. 

\end{defi}

\begin{rem}
    Being a torsor under a vector bundle over the $p$-torsion-free scheme $S$, $\mathbf U$ is $p$-torsion-free as well.
\end{rem}

\begin{rem}
Let $X$ be a smooth projective scheme over $\Z_p$. The (absolute) frobenius of the $\F_p$-variety $\overline X$, seldom lifts to an endomorphism $\Frob_X: X \to X$. Indeed, the existence of such a lift puts heavy constraints on the geometry of $\overline X$; see the nice paper \cite{D2}. In the current setting, we are now going to show that, over $\overline {\mathbf U} $, the (divided) verschiebung $\ver_{\overline V}$  acquires a natural splitting.

\begin{mot}
    Over a $p$-torsion-free base, if (the mod $p$) frobenius lifts, then it splits.
\end{mot} Hopefully, this  is clarified by  the explicit computations  coming next.
\end{rem}

\begin{lem}\label{LemDefis}
    Consider the composite \[ c:V^{(1)} \xrightarrow{\Ver_V} \Sym^p(V) \xrightarrow{\mu} \Gamma^p(V). \] Then, $c$ vanishes modulo $p$. Thus, there exists a unique $\mathcal O_{\mathbf U}$-linear arrow \[ \phi_V:V^{(1)} \to \Gamma^p(V), \] such that $p \phi_V=c$.
\end{lem}
\begin{dem}

Recall the formula giving    $\Ver_{\overline V}$: \[ \overline  V^{(1)} \to \Sym^p(\overline  V),\] \[\overline  x^{(1)} \mapsto \overline x ^p.\] The first claim follows, since \[\mu( \overline x ^p)=[\overline x]_1^p= p![\overline x]_p=0 \in \Gamma^p (\overline V). \]
For the second one, just recall that $\mathbf U$ is $p$-torsion-free.
    \end{dem}
\begin{lem}\label{LemSplitfrob}
    The mod $p$ reduction $ \overline \phi_V:  \overline V^{(1)} \to \Gamma^p( \overline V) $  splits $-\frob_{\overline V}.$
\end{lem}

\begin{dem}
    The claim is Zariski-local, so that one may assume that $S=\Spec(R)$ and $V=R^d$ (for $R$ a $p$-torsion-free ring). Denote by $(e_1,\ldots,e_d)$ the canonical basis of $V$.  Clearly, the homomorphism 
    $$\fonction{ \Ver_V '}{V^{(1)}}{\Sym^p( V)}{e_i^{(1)}}{ e_i^p}$$
  is another lift of  $\ver_{\overline V}$, defined over $S$. Hence, there exists a homomorphism of vector bundles over $\mathbf U $, \[ \delta:  V^{(1)} \to \Sym^p( V),\] such that $\Ver_V=\Ver_V'+p \delta$. Define the homomorphism of $S$-vector bundles 
  $$\fonction{ \phi'}{V^{(1)}}{ \Gamma^p(V)}{e_i^{(1)}}{ [e_i]_p.}$$
 One readily checks, that $\overline \phi'$ splits $\frob_{\overline V}$, and that \[\mu \circ \Ver_V'=p! \phi'.\] One then computes \[ p\phi_V= c= \mu \circ \Ver_V= p!\phi'+p \mu \circ \delta. \] Upon dividing by $p$, one finds  \[ \phi_V= (p-1)!\phi'+ \mu \circ \delta. \] Reducing mod $p$ then gives\[ \overline \phi_V=-\overline \phi' + \overline \mu \circ  \overline \delta, \]  so that  \[ \frob_{\overline V} \circ \overline \phi_V=- \frob_{\overline V} \circ \overline \phi' + \frob_{\overline V} \circ \overline \mu \circ  \overline \delta=- \Id_{\overline V^{(1)}} +0 , \] because $\frob_{\overline V} \circ \overline \mu=0 $. The claim is proved. 
\end{dem}

\begin{defi}
Define a natural linear arrow (where $d=\dim(V)$) \[ \nat:  \Lambda^{d-1} (\Gamma^p(V)) \to \Gamma^p( \Lambda^{d-1}(V)) \]  by the formula \[ [x_1]_p \wedge\ldots \wedge [x_{d-1}]_p \longmapsto [x_1 \wedge\ldots \wedge x_{d-1}]_p. \]

\end{defi}

\begin{rem}
    A word of explanation is in order, to justify that $\nat$ is well-defined. By gluing, one can assume that $S=\Spec(R)$ is affine. The functorial association 
    $$\fonctionnoname{V^{d-1}}{ \Lambda^{d-1} (\Gamma^p(V))}{(x_1,\ldots,x_{d-1})}{[x_1 \wedge\ldots \wedge x_{d-1}]_p}$$
    defines a polynomial law of $R$-modules, multi-homogeneous of multi-degree $(p,p,\ldots,p)$, in the sense of \cite{Fe}. By the universal property of divided powers, it gives rise to a multi-linear arrow of $R$-modules \[M: \Gamma^p(V)^{d-1} \to  \Lambda^{d-1} (\Gamma^p(V)),\] given on pure symbols by the formula\[ ([x_1]_p,\ldots,[x_{d-1}]_p) \longmapsto  [x_1 \wedge\ldots \wedge x_{d-1}]_p.\] One readily checks, that the expression on the right is anti-symmetric, and  vanishes whenever two $x_i$'s are equal. Let us check that  $M$ is actually alternating. This is clear if $\Gamma^p(V)$ is generated by  pure symbols $[x]_p$. Such is the case if $R=\Q$, e.g. by  item (3) of Lemma \ref{LemGenSym}. To conclude, let us explain how to reduce to this case. That $M$ is alternating,  is Zariski-local on $S$. We may thus assume that $V=R^d$ is a trivial vector bundle. Then, the situation is base-changed from $\Spec(\Z)$, so that w.l.o.g. we can further assume $R=\Z$. The statement can then be checked upon applying $(. \otimes_{\Z} \Q)$,  because the source and target of $M$ are free $\Z$-modules.
\end{rem}

\begin{lem}\label{Lemnat}
One works here mod $p$.\\
    Let $a_1 \in \Sym^p(\overline V)$ and $b_2,\ldots, b_{d-1} \in \Gamma^p(\overline  V)$. Then, one has \[ \overline  \nat( \overline  \mu(a_1) \wedge b_2 \wedge \ldots  \wedge b_{d-1}) \in \Im( \Sym^p(\Lambda^{d-1}(\overline  V)) \xrightarrow{\overline  \mu }\Gamma^p(\Lambda^{d-1}(\overline  V)) ).\]
    In short: if one of the  $(d-1)$ variables at the source of $\overline  \nat$  belongs to $\Im(\overline  \mu)$, then the image belongs to $\Im(\overline  \mu)$ as well.
\end{lem}

\begin{dem}
    On $\overline S$,  $(p-1)!=-1$ is invertible, so that by item (2) of Lemma \ref{LemGenSym} (for $n=p$ and $M=\overline V$), $\Im(\overline   \mu) \subset \Gamma^p (\overline  V)$ is spanned by the elements  (for $x,y \in \overline V$) \[z:=[x+y]_p-[x]_p-[y]_p.\] One then computes, for $x_2, \ldots, x_{d-1} \in \overline  V$: \[\overline  \nat( z \wedge [x_2]_p \wedge \ldots  \wedge  [x_d]_p)\] \[= [(x+y) \wedge x_2 \wedge\ldots \wedge x_{d-1}]_p   - [x \wedge x_2 \wedge\ldots \wedge x_{d-1}]_p  - [y \wedge x_2 \wedge\ldots \wedge x_{d-1}]_p.\] This indeed belongs to $\Im(\overline  \mu)$  (apply item (2) again, to $M= \Lambda^{d-1}(\overline  V)$). \\One concludes because symbols $[.]_p$ generate $\Gamma^p (\overline  V)$, by  item (3) of  Lemma \ref{LemGenSym}.
\end{dem}

\begin{defi}\label{defipsiV}
    Consider the composite
    \[ V^{(1)\vee} \otimes \Det(V) ^{-p}= \Lambda^{d-1} (V^{(1)}) \xrightarrow{\lambda^{d-1}(\phi_V)}  \Lambda^{d-1} (\Gamma^p(V)) \xrightarrow{\nat} \Gamma^p( \Lambda^{d-1}(V))=\Gamma^p(V^\vee ) \otimes \Det(V) ^{-p}. \] Upon twisting by $\Det(V)^p$, denote the resulting  arrow by \[\phi_{V^\vee}: V^{(1)\vee} \to\Gamma^p(V^\vee ), \]and its  dual  by \[\psi_V: \Sym^p(V) \to V^{(1)}. \] 
\end{defi}

\begin{lem}\label{Lemfrobsplit}
    The mod $p$ reduction $ \overline \phi_{V^\vee}:   \overline V^{\vee (1)} \to\Gamma^p( \overline V^\vee ) $ splits  $\frob_{\overline V^\vee}$. \\ Equivalently, $ \overline \psi_V$ splits  $\ver_{\overline V}$.
\end{lem}
\begin{dem}
 As in the proof of Lemma \ref{LemSplitfrob}, from which we retain notation, one may assume $S=\Spec(R)$ and $V=R^d$. Denote by $(e_1^\vee,\ldots, e_d^\vee)$ the basis dual to the canonical basis $(e_1,\ldots,e_d)$. Recall that \[e_i^\vee \otimes (e_1 \wedge  \ldots \wedge e_d) =(-1)^{i-1} e_1 \wedge \dots \wedge \widehat e_i \wedge \ldots \wedge  e_d,\] via the natural isomorphism \[V^\vee \otimes \Det(V) \simeq \Lambda^{d-1} (V),\] and where notation $ \widehat x$ means that $x$ is omitted.  Using this, one computes \[ (\nat \circ \lambda^{d-1}(\phi'))(e_i^{\vee (1)})= \nat(\lambda^{d-1}(\phi')((-1)^{i-1} e_1^{(1)} \wedge \dots  \wedge \widehat {e_i^{(1)}} \wedge \ldots \wedge e_d^{(1)}))\] \[ =(-1)^{i-1} \nat([e_1]_p \wedge \ldots  \wedge \widehat {[e_i]_p }\wedge \ldots \wedge [e_d]_p) =(-1)^{i-1} [e_1 \wedge \dots  \wedge \widehat e_i \wedge \ldots \wedge e_d]_p=[e_i^\vee]_p.\] [In the middle  steps of this computation, there should be a twist by $(e_1 \wedge \ldots \wedge e_d)^{-p} \in \Det(V)^{-p}$, which we ignored for clarity.]\\  Since $\overline \phi_{V}=\overline \phi'+\overline \mu \circ \overline \delta$, Lemma \ref{Lemnat} shows that  the difference \[\overline \nat \circ \lambda^{d-1}( \overline \phi_{V})- \overline \nat \circ  \lambda^{d-1}(\overline \phi')\] takes values in $\Im(\overline \mu)$. Since $\frob_{\overline V^\vee}$ vanishes on $\Im(\overline \mu)$, this implies   \[\frob_{\overline V^\vee} \circ \overline \nat \circ \lambda^{d-1}( \overline \phi_V)= \frob_{\overline V^\vee} \circ \overline \nat \circ  \lambda^{d-1}(\overline \phi')=\Id_{\overline V^\vee},\] where the last equality is readily checked on the basis $(e_i^{\vee (1)}),$ using the (mod $p$ reduction of the) computation above.
\end{dem}

\begin{lem}\label{LemDefiVerDual}
  Denote by $\Ver_{V^\vee}$ the composite (of $\mathcal O_{\mathbf U}$-linear maps) \[ V^{\vee (1)} \xrightarrow{ \phi_{V^\vee}} \Gamma^p(V^\vee) \xrightarrow{[x]_p \mapsto x^p} \Sym^p(V^\vee). \]  It is a lift of $ \ver_{ \overline V^\vee}$.
\end{lem}
\begin{dem}

Use Lemma \ref{Lemfrobsplit}, and the (mod $p$) factorisation of $[x]_p \mapsto x^p$, as  \[  \Gamma^p(\overline V^\vee) \xrightarrow{\frob_{\overline V^\vee}}  \overline V^{(1)}   \xrightarrow{\ver_{\overline V^\vee}}\Sym^p(\overline V^\vee).\]

\end{dem}

\begin{rem}\label{remdual}
    Altogether, over $\mathbf U=\mathbf U(V_2)$, there exist a natural linear map \[\psi_V: \Sym^p(V) \to V^{(1)},\] whose mod $p$ reduction splits $\ver_{\overline V}$. Also, there is   a natural linear map \[\phi_V: V^{(1)} \to \Gamma^p(V) ,\] whose mod $p$ reduction splits $\frob_{\overline V}$, and we have  seen that $\ver_{\overline V^\vee}$ lifts, to a natural $\mathcal O_{\mathbf U}$-linear map \[V^{\vee (1)} \xrightarrow{\Ver_{V^\vee}} \Sym^p(V^\vee).\] By universal property of the uplifting scheme, this lift is obtained from the universal one, by base-change via a natural $S$-morphism  \[\mathbf U(V_2) \to \mathbf U(V_2^\vee). \]
    Consequently, with regards to the existence of splittings and liftings, the situation here is `self-dual'. This is to be understood  in a weak sense, as the morphism above is not an isomorphism.
\end{rem}

The next Lemma states that, over $\mathbf U$, the $\W_2$-bundle $V_2$ arises as a push-forward of $\W_2(V)$- the universal lift of $V$ to a $\W_2$-Module, see  section \ref{secTeichLiftVB}.
\begin{lem}\label{lempushW2}
Recall the arrow $\psi_V$ of Definition \ref{defipsiV}.
Over $\mathbf U$, there is  a unique $\W_2$-linear arrow \[\Psi_V: \W_2(V) \to V_2,\] fitting into the  push-forward diagram (of extensions of  $\W_2$-Modules over $\mathbf U$)
     \[\xymatrix{\mathcal R \W_2(V): 0 \ar[r]  & \Frob_*(\Sym^p(V))\ar[r]\ar[d]^{\Frob_*(\psi_V)} &  \W_2(V)\ar[r] \ar[d]^{\Psi_V} &  V \ar[r] \ar@{=}[d]&  0\\ \mathcal R V_2: 0 \ar[r]  & \Frob_*(V^{(1)})\ar[r] &  V_2 \ar[r]&  V \ar[r] &  0.}\]
  
\end{lem}

\begin{dem}
Define the pushed-forward extension 
 \[\xymatrix{\mathcal E:= \Frob_*(\psi_V)_*(\mathcal R \W_2(V)): 0 \ar[r]  & \Frob_*(V^{(1)})\ar[r] &  \ast  \ar[r]&  V \ar[r] &  0.}\] The statement amounts to the existence of a \textit{unique} section of the Baer difference \[\xymatrix{\Delta:=(\mathcal E -\mathcal R V_2): 0  \ar[r]  & \Frob_*(V^{(1)})\ar[r] & D \ar[r]^\delta &  V \ar[r] &  0.}\]
Observe that, if such a section exists, it is unique  because, by Lemma \ref{LemExtW2Split}, \[ \Hom_{\W_2(\mathcal O_{\mathbf U})}( V,  \Frob_*(V^{(1)}))=0.\] As a consequence, the sought-for splitting may be constructed Zariksi-locally, allowing us to assume that $S$ (hence $\mathbf U$) is affine. Let us examine the situation mod $p$, i.e. upon base-change to $\overline  {\mathbf U}$. Consider the  push-forward diagram (of extensions of  $\W_2$-Modules over $\overline  {\mathbf U}$)
     \[\xymatrix{\mathcal R \W_2(\overline V): 0 \ar[r]  & \frob_*(\Sym^p( \overline  V))\ar[r]\ar[d]^{\frob_*(\overline  {\psi_V})} &  \W_2(\overline V)\ar[r] \ar[d] &  \overline V \ar[r] \ar@{=}[d]&  0\\ \overline{ \mathcal E}: 0 \ar[r]  & \frob_*(\overline V^{(1)})\ar[r] & \ast  \ar[r]&  \overline  V \ar[r] &  0.}\]
      Recall (see \cite{DCFL}, section 3.5) the connecting arrows \[ \kappa_{\mathcal R \W_2(\overline V)}:  \overline V  \to  \frob_*(\Sym^p( \overline  V)),\] resp. \[ \kappa_{\overline {\mathcal E}}:  \overline V  \to\frob_*(\overline V^{(1)})\] induced by the upper, resp. lower line of this diagram. By Lemma \ref{lemkappaad}, we know that $ \kappa_{\mathcal R \W_2(\overline V)}$ corresponds, via the ajdunction $(\frob^*,\frob_*)$, to the verschiebung \[\ver_{\overline V}: \overline V^{(1)} \to \Sym^p( \overline  V).\] Consequently, using compatibility of the formation of $\kappa$ to push-forwards,   \[ \kappa_{\overline {\mathcal E}}=(\frob_*(\overline  {\psi_V})\circ  \kappa_{\mathcal R \W_2(\overline V)})\] corresponds to $ (\overline  {\psi_V}\circ \ver_{\overline V})=\Id_{\overline V}$   (Lemma \ref{Lemfrobsplit}).  Thus, $\kappa_{\overline {\mathcal E}}$ is the  adjunction \[\mathrm{ad}_V:  \overline V  \to\frob_*(\overline V^{(1)})\] \[ x \mapsto x\otimes 1,\] which is also $\kappa_{\mathcal R\overline V_2}$. Using compatibility of $\kappa$ to Baer sum, it follows that \[\kappa_{\overline \Delta}=\kappa_{\overline {\mathcal E}}-\kappa_{\mathcal R \overline V_2}=0.\] This implies that $\overline \Delta$ is actually an extension of \textit{vector bundles} over  the affine scheme $\overline  {\mathbf U}$ (in other words, $\overline D$ is a vector bundle). A such, it splits (non-canonically). Let $\tilde g: \overline V \to \overline D $ be a splitting of $\overline {\Delta}$. Since $\mathbf U$ is affine, there exists an $\mathcal O_{\mathbf U}$-linear map \[ g: V \to D,\] such that $\overline g=\tilde g.$ From the assumptions, and the fact that $p$ is nowhere a zero-divisor on $\mathbf U$, one sees that there exists a unique endomorphism $\epsilon: V \to V,$ such that \[\delta \circ g=\Id_V+p \epsilon.\]  By Lemma \ref{LemExtW2Split},  $p \Delta$ splits, so that $p \epsilon$  factors through $\delta$. Consequently, so does $\Id_V$. In other words: $\Delta$ splits, as was to be shown. 
\end{dem}

\begin{lem}
    One may perform the  construction `dual' to that of Definition \ref{lempushW2}, using $\phi_V^\vee$ in place of $\psi_V$: over  $\mathbf U$, there is a natural  push-forward diagram 
     \[\xymatrix{\mathcal R \W_2(V): 0 \ar[r]  & \Frob_*(\Sym^p(V^\vee))\ar[r]\ar[d]^{\Frob_*(\phi_V^\vee)} &  \W_2(V^\vee)\ar[r] \ar[d] &  V^\vee \ar[r] \ar@{=}[d]&  0\\ \mathcal R V_2^\vee: 0 \ar[r]  & \Frob_*(V^{\vee(1)})\ar[r] &  V_2^\vee\ar[r]&  V^\vee \ar[r] &  0.}\]
     
\end{lem}

\begin{dem}
Same as that of Lemma  \ref{lempushW2}.
\end{dem}

\subsection{Uplifting w.r.t.  $\RR$.\\}

Let $r \geq 0$ be an integer. Let $q:=p^{r+1}$. Recall the ring scheme $\RR:=\W_2^{[r]}$ (Definition \ref{DefiW2r}), coming with a natural surjection of ring schemes over $\Z$, \[\pi: \W_{r+2} \to \RR.\]
Let $V$ be a vector bundle over the $p$-torsion-free scheme $S$, of constant rank $d\geq 1$. The goal here is to achieve the same contructions as in the previous section, using the ring scheme $\RR$ in place of $\W_2$. To do so, the (only, but significant) technical difficulty, is to generate a lift of 
$$\fonction{ \Ver_{\overline V}^r}{\overline V^{(r+1)}}{\Sym^q(\overline V)}{ x \otimes 1}{x^q.}$$
Indeed, the existence of a lift for $r=0$ (i.e. for $q=p$) does not at all imply that of a lift for  $r \geq 1$. To resolve this difficulty, \begin{center}
    \textit{ assume given a lift of $V$, to a $\W_{r+2}$-bundle $V_{r+2}$ over $S$.}
\end{center}

\begin{defi}

     Let $i\in \{1,\ldots,r+1\}$  be an integer.   As usual, denote by \[V_i:=(\rho^{r+2-i})^*(V_{r+2})\] the reduction of $V_{r+2}$, to a $\W_i$-bundle over $S$.\\
     For $i\in \{0,\ldots,r\}$, define the $\W_2$-bundle \[V_2^{(r-i)}:=(\Frob^{r-i})^*(V_{r+2-i}),\] and the uplifting scheme of the $\W_2$-bundle $\Sym^{p^i}(V_2^{(r-i)})$, \[ \mathbf U_i := \mathbf U(\Sym^{p^i}(V_2^{(r-i)})) \to S. \]
     Define the uplifting scheme of $V_{r+2}$,  as \[ \mathbf U=\mathbf U(V_{r+2}):=\prod_{i=0}^{r} \mathbf U_i  \to S\](product fibered over $S$.) 
\end{defi}
\begin{lem}\label{defipsiR}
  Let $i\in \{0,\ldots,r\}$  be an integer. \begin{enumerate}
      \item{Over $\mathbf U_i$, the  (divided) verschiebung
      $$\fonction{\ver}{ \Sym^{p^i}(\overline V^{(r+1-i)})}{\Sym^{p^{i+1}}(\overline V^{(r-i)})}{ x \otimes 1}{ x^p}$$
   acquires a natural lift, to an $\mathcal O_{\mathbf U_i}$-linear map denoted by \[ \Ver_i: \Sym^{p^i}(V^{(r+1-i)}) \to \Sym^{p^{i+1}}(V^{(r-i)}).\] } \item{ There exists  a natural $\mathcal O_{\mathbf U_i}$-linear map  \[ \psi_i: \Sym^{p^{i+1}}(V^{(r-i)}) \to \Sym^{p^i}(V^{(r+1-i)}),\]  such that, modulo $p$, \[(\overline \psi_i \circ \ver^{i+1})=\ver^i: \overline  V^{(r+1)} \to \Sym^{p^i}(\overline  V^{(r+1-i)}) .\]}  \item{Over $\mathbf U$, denote by \[\psi: \Sym^{p^{r+1}}(V) \to V^{(r+1)}\] the composite \[ \Sym^{p^{r+1}}(V) \xrightarrow{\psi_r}  \Sym^{p^{r}}(V^{(1)}) \xrightarrow{\psi_{r-1}}  \Sym^{p^{r-1}}(V^{(2)})  \xrightarrow{\psi_{r-2}}  \ldots  \xrightarrow{\psi_0}  V^{(r+1)} .\] Then, $\overline \psi$ is a retraction of $ \overline V^{(r+1)} \xrightarrow{\ver^{r+1}}\Sym^{p^{r+1}}(\overline V), $ over $\overline {\mathbf U}$. }
  \end{enumerate}

\end{lem}

\begin{dem}

Denote by \[\ver'_i:\Sym^{p^i}(\overline V^{(r+1-i)})  \to \Sym^p(\Sym^{p^i}(\overline V^{(r-i)}) ) \] the  verschiebung  of the vector bundle $\Sym^{p^i}(\overline V^{(r-i)}) $.  By definition of $\mathbf U_i$, it lifts over $\mathbf U_i$, to a linear map   \[ \Ver'_i: \Sym^{p^i}(V^{(r+1-i)}) \to \Sym^p(\Sym^{p^{i}}(V^{(r-i)})).\] 
Define $\Ver_i$ as the composite \[\Sym^{p^i}(V^{(r+1-i)}) \xrightarrow{\Ver'_i} \Sym^p(\Sym^{p^{i}}(V^{(r-i)})) \xrightarrow{\mu} \Sym^{p^{i+1}}(V^{(r-i)}).\]
One computes the  mod $p$ reduction of $\Ver_i$, as  \[\Sym^{p^i}(\overline V^{(r+1-i)}) \xrightarrow{\ver'_i} \Sym^p(\Sym^{p^{i}}(\overline V^{(r-i)})) \xrightarrow{\overline \mu} \Sym^{p^{i+1}}(\overline V^{(r-i)}).\] By (the dual formulation of) Lemma \ref{lemnatfrob}, this composite is none other than $\ver$, the verschiebung of $\overline V$. This proves the first claim. For item (2), denote by \[ \psi_i': \Sym^p(\Sym^{p^{i}}(V^{(r-i)})) \to  \Sym^{p^i}(V^{(r+1-i)}) \] the linear map introduced in Definition \ref{defipsiV}
 (applied to the $\W_2$-bundle  $\Sym^{p^i}(V_2^{(r-i)})$). By Lemma \ref{Lemfrobsplit}, its mod $p$ reduction $\overline \psi_i'$ splits $ \ver_i'$. Define $\psi_i$ as the composite \[\Sym^{p^{i+1}}(V^{(r-i)}) \xrightarrow{\sigma^p} \Sym^p(\Sym^{p^{i}}(V^{(r-i)})) \xrightarrow{ \psi_i'} \Sym^{p^i}(V^{(r+1-i)})\] (see Definition \ref{defisigmasym} for  $\sigma^p$.)

 Modulo $p$, using  item (2) of Lemma \ref{lemnatfrob}, one gets \[\overline \psi_i \circ \ver^{i+1}= \overline \psi_i' \circ \overline \sigma^p \circ \ver^{i+1}= \overline \psi_i' \circ\ver_i' \circ \ver^i=\ver^i,\]  proving the claim.  Checking (3) is straightforward: using (2) for all $i$, one computes that the composite $(\overline \psi \circ \ver^{r+1}):\overline V^{(r+1)} \to \overline V^{(r+1)}$equals  \[  \overline \psi_0 \circ  \ldots \circ  \overline \psi_r \circ \ver^{r+1}=\overline \psi_0 \circ  \ldots \circ  \overline \psi_{r-1} \circ \ver^r=\overline \psi_0 \circ  \ldots \circ  \overline \psi_{r-2} \circ \ver^{r-1}= \ldots = \psi_0 \circ \ver=\Id.\]

\end{dem}

item (3) of the preceding Lemma makes possible the following Definition. The existence of the push-forward diagram there, is obtained in the same way as in the proof of Lemma \ref{lempushW2}, \textit{mutatis mutandis}. We omit the straightforward details.
\begin{defi}\label{defiUS}
Let $S$ be a scheme where $p$ is nowhere a zero-divisor. Let $r \geq 0$ be an integer, and let $V_{r+2}$ be a $\W_{r+2}$-bundle over $S$. Set $q:=p^{r+1}$, and denote by \[\pi: \W_{r+2} \to \RR\] the natural surjection of ring schemes over $\Z$ (see Definition \ref{DefiW2r}). Define \[V_{\RR}:=V_{r+2} \otimes_\pi \RR.\] It is an $\RR$-bundle over $S$, lifting the vector bundle $V:=\rho(V_{r+2})$. Observe that $V^{[1]}=V^{(r+1)}$ (Lemma \ref{Lembrpr}.)\\
   Over the uplifting scheme \[\mathbf U=\mathbf U(V_{r+2}) \to S,\]  there is  a unique $\RR$-linear arrow \[\Psi: \RR (V) \to V_{\RR}, \]fitting into the  push-forward diagram (of extensions of  $\RR$-Modules)
     \[\xymatrix{\mathcal R \RR(V): 0 \ar[r]  & (\Frob_\RR)_*(\Sym^q(V))\ar[r]\ar[d]^{\Frob_*(\psi)} &  \RR(V)\ar[r] \ar[d]^\Psi &  V \ar[r] \ar@{=}[d]&  0\\ \mathcal R V_{\RR}: 0 \ar[r]  & (\Frob_\RR)_*(V^{[1]})\ar[r] &  V_{\RR} \ar[r]&  V \ar[r] &  0,}\]

     where $\psi$ is defined in item (3) of Lemma \ref{defipsiR}.
\end{defi}

 \subsection{Uplifting for duals.\\}\label{secdualUL}

Keep notation and assumptions of Definition \ref{defiUS}.\\
Recall Remark \ref{remdual}, that  expresses  a `self-dual' property of uplifing schemes  w.r.t.  $\W_2$. It is given by a natural morphism of $S$-schemes $\mathbf U(V_2) \to \mathbf U(V^\vee_2)$.\\ Applying this to  all $\mathbf U_i$'s above,  thus yields a natural morphism of $S$-schemes $ \mathbf U(V_{r+2}) \to \mathbf U(V_{r+2}^\vee)$. Consider Definition \ref{defiUS}, applied to $V_{r+2}^\vee$. Upon base-change via this morphism,  one gets the following dual counterpart.

\begin{it}
 \noindent Over $\mathbf U=\mathbf U(V_{r+2})$,   there is a natural linear map \[\psi': \Sym^q(V^\vee)  \to V^{\vee[1]},\] whose mod $p$ reduction is a retraction of $\ver^{r+1}:  \overline V^{\vee (r+1)} \to  \Sym^q(\overline V^\vee)$. \\Moreover, there is a unique $\RR$-linear arrow \[\Psi': \RR (V^\vee) \to V^\vee_{\RR}, \]fitting into the  push-forward diagram (of extensions of  $\RR$-Modules over $\mathbf U$)
     \[\xymatrix{\mathcal R \RR(V^\vee): 0 \ar[r]  & (\Frob_\RR)_*(\Sym^q(V^\vee))\ar[r]\ar[d]^{\Frob_*(\psi')} &  \RR(V^\vee)\ar[r] \ar[d]^{\Psi'} &  V^\vee \ar[r] \ar@{=}[d]&  0\\ \mathcal R V^\vee_{\RR}: 0 \ar[r]  & (\Frob_\RR)_*(V^{\vee[1]})\ar[r] &  V^\vee_{\RR} \ar[r]&  V^\vee \ar[r] &  0.}\]

\end{it}  


\subsection{Uplifting for tensor products.\\}\label{sectensUL}Keep notation and assumptions of Definition \ref{defiUS}.\\ Let $V'_{r+2}$ be another $\W_{r+2}$-bundle over $S$.  Denote its uplifting scheme by \[\mathbf U'=\mathbf U(V'_{r+2}) \to S.\] Set \[V_{r+2}'':=V_{r+2} \otimes_{\W_{r+2}} V_{r+2}', \] and define the fibered product \[\mathbf U'':=\mathbf U \times_S \mathbf U' \to S.\]
Over $\mathbf U''$, there are the arrows built above, \[\psi: \Sym^q (V) \to V^{[1]}\]and \[\psi':  \Sym^q (V') \to V^{'[1]}.\]  Define an $\mathcal O_{\mathbf U''}$-linear map $\psi''$, as the composite \[\Sym^q (V'') \xrightarrow{\nat} \Sym^q (V )\otimes \Sym^q (V') \xrightarrow{\psi \otimes \psi'} V^{[1]}\otimes V^{'[1]}=V^{''[1]}, \] where the natural quotient $\nat$ is given by the formula \[ (x_1 \otimes x'_1)(x_2 \otimes x'_2) \ldots (x_q\otimes x'_q) \longmapsto  (x_1 x_2 \ldots x_q) \otimes  (x'_1 x'_2 \ldots x'_q) .\]

 One readily checks that, over $\overline {\mathbf U''}$, $\overline {\psi''}$ is a retraction of \[\overline {V^{''}}^{(r+1)} \xrightarrow{\ver^{r+1}}\Sym^q(\overline { V''}). \]

Consequently, one can perform, over $\mathbf U''$, the same construction as in Definition \ref{defiUS}: there is  a unique $\RR$-linear arrow $\Psi''$, fitting into the  push-forward diagram 
     \[\xymatrix{\mathcal R \RR(V''): 0 \ar[r]  & (\Frob_\RR)_*(\Sym^q(V''))\ar[r]\ar[d]^{\Frob_*(\psi'')} &  \RR(V'')\ar[r] \ar[d]^{\Psi''} &  V'' \ar[r] \ar@{=}[d]&  0\\ \mathcal R V''_{\RR}: 0 \ar[r]  & (\Frob_\RR)_*(V^{''[1]})\ar[r] &  V''_{\RR} \ar[r]&  V'' \ar[r] &  0.}\]

     \noindent  This process clearly extends, to the case of tensor products of three or more bundles.
   
\section{Suspension.}\label{SecSuspension}

\noindent We describe a process for lifting Hochschild $1$-cocycles, which we call suspension. It is a universal construction,  performed via  \textit{group-change}, rather than  base-change.  It  is a major ingredient to build  the Uplifting Pattern. 
\begin{rem}
    It would be interesting to study suspension, in the more elementary context of cohomology of discrete $G$-modules, where $G$ is a (pro)finite group.
\end{rem}

  \subsection{Suspension of a $1$-cocycle.}\label{DefiSusp}\hfill \\
Let $S$ be an affine scheme. Let $\Gamma \to S$ be an affine  flat $S$-group scheme.
  Let $$(E): 0 \to U \to   \tilde M  \xrightarrow{\pi}   M\to 0$$ be an admissible extension of $(\Gamma,S)$-groups. Assume that $U$ is an affine  flat $S$-group scheme.\\
 Consider a $1$-cocycle over $S$, $$ C \in Z^1(\Gamma,M).$$  We would like to produce a change of affine flat $S$-group scheme  \[\gamma: \tilde \Gamma \to \Gamma, \] which is universal for the property, that  the cocycle $$ \gamma^*(C) \in Z^1( \tilde \Gamma ,M)$$ lifts via $\pi_*$, to a cocycle  $$\tilde C \in Z^1(\tilde \Gamma,\tilde M).$$
  
  This is achieved by a natural (actually, obvious) construction, that does not rely on connecting maps nor on $2$-cocycles, but rather  on the
  
\begin{mot}
Fiber product of a group homomorphism and a $1$-cocycle is a group.
\end{mot}
To understand this, form the fibered product of $S$-functors
\[ \xymatrix{ \tilde \Gamma \ar[r]^-{\gamma} \ar[d] & \Gamma  \ar[d]^{C} \\  \tilde M \ar[r]  &   M.}\] 

Endow $\tilde \Gamma$ with the structure of an  $S$-group, point-wise given by  \[ (\gamma,x)(\gamma',x'):=(\gamma \gamma',x+\gamma(x')),\] for $\gamma,\gamma' \in \Gamma$ and $x,x' \in \tilde M$. It is readily checked,  that this formula indeed defines a group law, and that $ \gamma$ has the desired lifting property.

There is a pull-back  diagram, with admissible rows  \[\xymatrix{(E \tilde \Gamma): 0 \ar[r] &  U   \ar[r] \ar@{=}[d] &\tilde \Gamma  \ar[r]^-\gamma  \ar[d]^-{\tilde  C}  & \Gamma  \ar[r]  \ar[d]^-C & 1 \\  (E): 0 \ar[r] &  U  \ar[r] &  \tilde M  \ar[r] &   M \ar[r] & 1. }\]
Since $\Gamma$ and $U$ are affine flat  $S$-schemes, so is $\tilde \Gamma$.
\begin{rem}
Later on, in typical applications,  $\Gamma$ is smooth over $S$, and $(E)$ is an extension of  smooth affine $S$-group schemes, with $U$ a vector group (i.e. the affine space of a vector bundle over $S$). In that situation,  the diagram above shows that $(E \tilde \Gamma)$ is an extension of  smooth affine $S$-group schemes, as well.
\end{rem}
  
  \begin{rem}(The commutative case.)\\\label{remcomm}
\noindent Keep notation above. Assume that $\tilde M$ is commutative. Denote by $$ [C] \in H^1(\Gamma,M)$$ the class of $C$. By Proposition \ref{PropFactSyst}, the admissible extension $E \tilde \Gamma$ has a class \[ [E \tilde \Gamma] \in H^2(\Gamma, U).\] From the definition of factor systems, one readily checks that $$[E \tilde \Gamma]=\beta_E ([C]),$$ where $$ \beta_E: H^1(\Gamma,M) \to H^2(\Gamma,U) $$ is the connecting map induced by $(E)$ in cohomology.
\end{rem}
 
 \subsection{Application: lifting extensions.}\label{Appsusp}

     Lifts are considered  here w.r.t. the ring scheme starring in this text, namely $\RR=\W_2^{[r]}$, for some $r \geq 0$. The same constructions work for other ring schemes.

 Let $S$ be a scheme.  Let $\Gamma$ be an affine flat $S$-group scheme. Consider   a  $\Gamma$-scheme $$f:X\to S.$$[Typically, in \cite{DCF3}, $S=\Spec(\Z)$, $\Gamma=\GL_n$, and $f$ is a smooth morphism.]\\  Let 
$$ \mathcal E: 0 \to   V\to E \to  \mathcal O_X\to 0 $$
be a geometrically split extension of $\Gamma$-bundles over $X$. By Proposition \ref{HochGT}, upon choosing a geometric splitting, (the iso class of) $\mathcal E$ is determined by a Hochschild $1$-cocycle over $S$,
 $$C: \Gamma \to\Pi_f(\A_X(V)).$$ 
 
 \begin{defi}
Let $V_{r+2}$ be a lift of $V$, to a $\Gamma \W_{r+2}$-bundle over $X$.  Denote by \[\pi: \W_{r+2} \to \RR\] the natural surjection (of ring schemes over $\Z$) and set \[ V_\RR:=V_{r+2} \otimes_{\pi} \RR;\] it is a lift of $V$, to a $\Gamma \RR$-bundle over $X$.     
 \end{defi}

\begin{defi}\label{defipiadmi}

Consider the reduction sequence of the $\Gamma \RR$-bundle  $V_\RR$, over $X$:  \[ \mathcal R: 0 \to (\Frob_\RR)_*( V^{[1]})  \to  V_\RR \to  V \to 0. \]
   We assume that the sequence of $S$-group functors  \[\Pi_f(\mathcal R):  0 \to \Pi_f(\A_X( V^{[1]}))  \to \Pi_f(\A_X(V_\RR ))\to \Pi_f(\A_X(V)) \to 0 \]

   is an admissible exact sequence of affine flat $S$-group schemes.
 \end{defi}
      Consider lifting $\mathcal E$, to a geometrically split  extension of $\Gamma \RR$-bundles $$ \mathcal E_\RR: 0 \to   V_\RR\to E_\RR \to \RR (\mathcal O_X) \to 0, $$ whose middle term $E_\RR$ is not prescribed.

 \begin{defi}\label{DefiLiftC}

     Apply the suspension process to the $1$-cocycle $C$, w.r.t. the  extension $\Pi_f(\mathcal R)$ of Definition \ref{defipiadmi}.  One gets an extension of  affine flat $S$-group schemes
 \[ 0 \to \Pi_f(\A_X(V^{[1]}))  \to \Gamma' \stackrel {\gamma} \to \Gamma \to 1.\] It is such that  the cocycle $$\gamma^*(C):  \Gamma' \to  \Pi_f(\A_X(V))$$ admits a tautological lift, to a cocycle  $$C_\RR:\Gamma' \to  \Pi_f(\A_X(V_\RR)),$$ which determines a    geometrically trivial extension of $\Gamma' \RR$-bundles lifting   $\gamma^*(\mathcal E)$, $$ \mathcal E_\RR: 0 \to   V_\RR\to E_\RR \to   \RR (\mathcal O_X) \to 0. $$ It is  universal in the following sense. For every affine morphism $T \xrightarrow{t} S$, and for every morphism of  $T$-group functors $H \xrightarrow{h} \Gamma_T$, there is a natural bijection between  \begin{enumerate}
     \item{The set of isomorphism classes of lifts of $h^*(t^*(\mathcal E))$, to an extension of $H \RR$-bundles  over $T$, reading as  $ \mathcal E_\RR$ above.} \item{The set of homomorphism of $T$-group functors $ H\to \Gamma'_T$, lifting $h$ through $ \Gamma'_T \xrightarrow{\gamma_T} \Gamma_T$, modulo conjugation by $ \Pi_f(\A_X( V^{[1]})) (T)$. }
 \end{enumerate} 
 \end{defi}

\section{The Uplifting  Pattern.}\label{secUP}
Suspension and uplifting interact well, combining to  produce a useful lifting result: The Uplifting Pattern. To simplify its exposition, we first ignore cyclotomic twists. The general case  will be taken care of in section \ref{secnontrivcyclo}.\\ Let $S$ be an \textit{affine} scheme, let $\Gamma$ be a smooth affine $S$-group scheme, and let 
$$ \mathcal E: 0 \to   V\to E \to  \mathcal O_S\to 0 $$
be a (geometrically split) extension of $\Gamma$-bundles over $S$,  determined (up to isomorphism) by a Hochschild $1$-cocycle over $S$,
 $$C: \Gamma \to \A_S(V).$$ 

Recall the reduction sequence  (where $q:=p^{r+1}$) \[\mathcal R\RR( V):  0 \to \A_S(\Sym^q(V))  \to \A_S( \RR( V))\xrightarrow{\rho_V} \A_S(V) \to 0 ,\] considered here as an extension  of smooth affine $S$-group schemes.

 \begin{defi}\label{DefiLiftRV}

Consider lifting the extension $\mathcal E$, via the arrow $\rho_V$ above, to an  extension of $\Gamma \RR$-Modules over $S$, of the shape $$ \RR(\mathcal E): 0 \to   \RR( V)\to \ast \to \RR (\mathcal O_S) \to 0. $$ 

Proceeding as in Definition \ref{DefiLiftC},   one gets a universal (admissible) extension of  affine smooth $S$-group schemes
 \[ 0 \to\A_S(\Sym^q(V))  \to \Gamma' \stackrel {\gamma} \to \Gamma \to 1,\]  together with an extension of $\Gamma' \RR$-Modules lifting   $\gamma^*(\mathcal E)$, reading as $ \RR(\mathcal E)$ above.
 \end{defi}

\begin{rem}
Denote by $\P(V) \xrightarrow{f} S$ the projective bundle of $V$. 
    From the construction of $\RR( V)$ in section \ref{secTeichLiftVB}, the lifting problem in definition \ref{DefiLiftRV}, is equivalent to the following one, over $\P(V)$.\\ Form the push-forward diagram (of extensions of $\Gamma$-bundles over $\P(V)$)
    
    \[\xymatrix {b^*(\mathcal E) :0 \ar[r] &b^*(V) \ar[r] \ar[d]^\nat  & b^*(E) \ar[r] \ar[d]& \mathcal O_{\P(V)} \ar[r] \ar@{=}[d] & 0 \\  \mathcal P :0 \ar[r] & \mathcal O_{\P(V)}(1) \ar[r] & \ast \ar[r] & \mathcal O_{\P(V)} \ar[r] & 0 .  }\] Consider lifting $\mathcal P$, to an  extension of $\Gamma \RR$-bundles over $\P(V)$, of the shape \[ \RR(\mathcal P): 0 \to  \RR(\mathcal O_{\P(V)}(1))\to \ast \to \RR (\mathcal O_{\P(V)}) \to 0. \]
    Indeed, one passes from $ \RR(\mathcal P)$ to $ \RR(\mathcal E)$ by applying $\Pi_f(.)$;  details are as in the proof  of Proposition \ref{proplifthom}.

\end{rem}
\begin{defi}

From now on, assume that $p$ is nowhere a zero-divisor on $S.$\\
 \noindent Let $V_{r+2}$ be a $\Gamma \W_{r+2}$-bundle over $S$, lifting the  $\Gamma$-bundle $V$.  Set \[ V_\RR:=V_{r+2} \otimes_{\pi} \RR;\] it is a lift of $V$, to a $\Gamma \RR$-bundle over $S$.\\
    Denote by \[\U=\U(V_{r+2} ) \xrightarrow{u} S\] its uplifting scheme; see Definition \ref{defiUS}. 
\end{defi}

\begin{prop}(The Uplifting Pattern, for the trivial cyclotomic twist.)\label{propULP}\\
    Upon group-change via $\gamma$, and base-change via  $u$, the  extension $(\mathcal E)$ lifts, to an extension of $\Gamma' \RR$-bundles over $\U$, reading as $$ \mathcal E_\RR: 0 \to   V_\RR\to \ast \to   \RR (\mathcal O_{\mathbf U}) \to 0. $$ 
\end{prop}

\begin{dem}
   Recall Definition \ref{defiUS}- especially the arrow $\Psi$. The pushed-forward extension $ \mathcal E_\RR:=\Psi_*(u^*( \RR(\mathcal E)))$ does the job.
\end{dem}
\subsection{Relation to cyclotomic pairs.\\}
Recall that $p$ is nowhere a zero-divisor on $S.$\\
Here we work upon reduction to $\overline S=S \times_{\Spec(\Z)} \Spec(\F_p)$.
 \begin{prop}\label{proplifthom}
   Let $G$ be a profinite group, such that the pair $(G,\Z/p^2)$ is $(1,1)$-cyclotomic. Consider a continuous (=whose kernel is open) group homomorphism \[\sigma: G \to \Gamma( S).\]  Then, for $r$ large enough, the mod $p$ reduction \[\overline \sigma: G \to \overline \Gamma( \overline S)\] factors through  the homomorphism $\gamma$ of Definition \ref{DefiLiftRV}. In other words, for $r >>0$, there exists  a continuous homomorphism $\sigma'$, such that $\overline \sigma$ equals the composite    \[ G  \xrightarrow{\sigma'} \Gamma'(\overline S)\xrightarrow{\gamma(\overline S)} \Gamma(\overline S).\]
 \end{prop}

\begin{dem}
Recall the definition of $\Gamma'$. Let $B:= \P(V) \xrightarrow{f} S$ be the projective bundle of $V$. Over $ B$, form the push-forward diagram (of extensions of  $\Gamma$-bundles over $B$)

     \[\xymatrix{f^*(\mathcal E): 0 \ar[r]  & f^*(V)\ar[r]\ar[d]^{\nat} & f^*(E) \ar[r] \ar[d] & \mathcal O_B \ar[r] \ar@{=}[d]&  0\\ \mathcal P: 0 \ar[r]  & \mathcal O_B(1) \ar[r] & P \ar[r]&  \mathcal O_B \ar[r] &  0.}\]
[Here $\mathcal O_B(1)$ denotes  the twisting sheaf on $B=\P(V)$, and has nothing to do with the cyclotomic twist $\Z/p^2(1)$, which is trivial  by assumption.]
    The extension $\mathcal P$ thus defined is geometrically split (because so is $\mathcal E$, over the affine base $S$).  Via  the given homomorphism $\sigma$, consider it as a  geometrically split extension of $G$-bundles over $B$.  Proposition \ref{propcyclolift} then applies-  over the perfection $\overline B^\perf$. The conclusion,  over $\overline B$, reads like this.  For $r \geq 0$ large enough,  the $r$-th frobenius twist  \[\xymatrix{ \overline {\mathcal P}^{(r)}: 0 \ar[r]  & \mathcal O_{ \overline B}(p^r ) \ar[r] & \overline  P^{(r)} \ar[r]&  \mathcal O_{ \overline B} \ar[r] &  0}\] lifts, to a (geometrically trivial)  extension of $G \W_2$-bundles over $\overline B$, \[\xymatrix{ 0 \ar[r]  &\W_2( \mathcal O_{ \overline B} (p^r ) )\ar[r] & \ast \ar[r]& \W_2( \mathcal O_{ \overline B} ) \ar[r] &  0.}\] Fix such an $r$, and take  $\RR:=\W_2^{[r]}$ for that $r$. Using the description of  
    $\RR$ over $\F_p$ (see Lemma \ref{LemR2overFp}), one gets the following result. Over $ \overline B$, the extension $\overline {\mathcal P}$ lifts, to a (geometrically trivial) extension of $G \RR$-bundles  \[\xymatrix{ 0 \ar[r]  &\RR( \mathcal O_{\overline B} (1 ) )\ar[r] & \ast \ar[r]& \RR( \mathcal O_{ \overline B}) \ar[r] &  0.}\] Applying $\Pi_{\overline f_*}(.)$ to this extension produces, over $\overline S$,  a lift of the extension $\overline {\mathcal E}$, to an extension of $G \RR$-Modules   \[\xymatrix{ 0 \ar[r]  &\RR(\overline V )\ar[r] & \ast \ar[r]& \RR( \mathcal O_{ \overline S} ) \ar[r] &  0.}\] 

  Universal property of suspension, then yields the sought-for homomorphism $\sigma'$.
\end{dem}

\subsection{Implementing a cyclotomic twist.\\}\label{secnontrivcyclo}
Notation here is as in the beginning of section \ref{secUP}: $S$ is an affine scheme,  $\Gamma$ is a smooth affine $S$-group scheme, and $ \mathcal E: 0 \to   V\to E \to  \mathcal O_S\to 0 $ is an extension of $\Gamma$-bundles on $S$.\\To deal with a non-trivial cyclotomic twist, the first step is to define its  `scheme-theoretic' counterpart.
To do so, recall from Definition \ref{defiGam}, the $\Z$-group scheme \[\G_{a,m}=\Ker(\RR ^\times \xrightarrow{\rho }\W_1^\times).\]

\begin{defi}
 Over $\Spec(\Z)$, define a $\G_{a/m} \RR$-line bundle $\RR[1]$ as follows. \\As an $\RR$-line bundle,  $\RR[1]=\RR(\Z)$, and the non-trivial $\G_{a/m}$-action is the natural one, obtained by restricting the natural multiplicative action of $\RR^\times$.\\
 Observe that $\RR(\Z)[1]$ is a lift of the trivial $\G_{a/m}$-line bundle  $\Z$.\\
 Put \[\Z\{q\}:=\Frob_\RR^*(\RR[1]);\] it is a  $\G_{a/m}$-line bundle  over $\Z$.
 \end{defi}

\begin{rem}
 Recall that the $\F_p$-group scheme $\G_{a/m} \times_{\Z} \F_p$  is the additive group $\G_a$,  whereas the $\Z[\frac 1 p]$-group scheme $\G_{a/m} \times_{\Z} \Z [ \frac 1 p]$  is the multiplicative group $\G_m$. Accordingly,  $\Z\{q\} \otimes_{\Z} \F_p$, is the trivial $\G_a$-line bundle ($=\F_p$), whereas the $\G_m$-action on $\Z\{q\} \otimes_{\Z} \Z [ \frac 1 p]$ is non-trivial, given by the character $\G_m \xrightarrow{x \mapsto x^q} \G_m$.
\end{rem}
 \begin{defi}

     Put $\Gamma_1:= \Gamma \times_{\Spec(\Z)} \G_{a/m}$, considered as a smooth affine $S$-group scheme. Henceforth, $\Gamma$-equivariant structures are considered as $\Gamma_1$ ones, via the projection $\Gamma_1 \to \Gamma$. If $M$ is a $\Gamma_1 \RR$-Module (over a $\Gamma_1$-scheme $X \to S$), put \[M[1]:=M \otimes_{\RR(\Z)} \RR[1].\] Similarly, if  $V$ is a $\Gamma_1$-bundle, put \[V\{q\}:=V \otimes_{\Z} \Z\{q\}.\] In both cases, the $\Gamma_1$-action on the second factor of $(.\otimes.)$ is given by the projection $\Gamma_1 \to \G_{a/m}$.
 \end{defi}

 One can now formulate a twisted version of Definition \ref{DefiLiftRV}:
\begin{defi}\label{DefiLiftRV1}
 Denote by $\mathcal E_1$ the extension $\mathcal E$, viewed as an extension of $\Gamma_1$-bundles.
Consider lifting $\mathcal E_1$, via the arrow \[ \rho_V[1]: \RR(V)[1] \to V[1]=V,\] to an  extension of $\Gamma_1 \RR$-Modules over $S$, of the shape $$ \RR(\mathcal E_1): 0 \to   \RR( V)[1]\to \ast \to \RR (\mathcal O_S) \to 0. $$ 

\noindent  Taking into account the extra presence of twists,  one then proceeds as in Definition \ref{DefiLiftC}. As a result, there is a  universal (admissible) extension of  affine smooth $S$-group schemes
 \[ 0 \to\A_S(\Sym^q(V)\{q\})  \to \Gamma_1' \stackrel {\gamma_1} \to \Gamma_1\to 1,\]  together with an extension of $\Gamma_1' \RR$-Modules lifting   $\gamma_1^*(\mathcal E_1)$, reading as $ \RR(\mathcal E_1)$ above.
 \end{defi}
Likewise, one gets:
\begin{prop}(The Uplifting Pattern.)\label{propULP1}\\
    Upon group-change via $\gamma_1$, and base-change via  $u$, the  extension $(\mathcal E_1)$ lifts, to an extension of $\Gamma_1' \RR$-bundles over $\U$, reading as $$ \mathcal E_{1,\RR}: 0 \to   V_\RR [1]\to \ast \to   \RR (\mathcal O_{\mathbf U}) \to 0. $$ 
\end{prop}
\begin{dem}
    Same as Proposition \ref{propULP}, using the arrow $\Psi[1]: \RR(V)[1] \to V_\RR[1] $.
\end{dem}

   Finally, let $(G,\Z/p^2(1))$ be a $(1,1)$-cyclotomic pair. Define \[\Z/p^2[1]:=\Z/p^2(1) \otimes_\Z \W_2(\F_p(-1)).\]
   
   \begin{lem}
   The pair  $(G,\Z/p^2[1])$ is  $(1,1)$-cyclotomic, and $\F_p[1]=\F_p$.     
   \end{lem}
   
   \begin{dem}
    Via  the classical restriction/corestriction argument (as in \cite{DCF1}, Remark 6.4), this can be checked upon replacing $G$ by  the kernel $G_0 \subset G$ of the $G$-action on $\F_p(1)$, because the index $[G:G_0]$ is prime-to-$p$. But then  $\Z/p^2[1]=\Z/p^2(1)$, so that the claim is obvious.   
   \end{dem}
   
 Using $\Z/p^2[1]$ in place of $\Z/p^2(1)$, Proposition \ref{proplifthom} generalises as follows.
 \begin{prop}\label{proplifthom1}
 Consider a continuous homomorphism \[\sigma: G \to \Gamma( S).\]  Denote by \[ \chi: G \to  (\Z/p,+)=\G_{a/m}(\F_p) \subset \G_{a/m}(\overline S) \] the character giving the action of $G$ on $\Z/p^2[1]$. Set \[\sigma_1:=(\overline \sigma,\chi): G \to \Gamma_1(\overline  S).\]    Then, for $r$ large enough, $\sigma_1$ factors through  the homomorphism $\gamma_1$ of Definition \ref{DefiLiftRV1}. In other words, for $r >>0$, there exists  a continuous homomorphism $\sigma_1'$, such that $\sigma_1$ equals the composite    \[ G  \xrightarrow{\sigma_1'} \Gamma_1'(\overline S)\xrightarrow{\gamma_1(\overline S)} \Gamma_1(\overline S).\]
 \end{prop}

\begin{dem}
 Observe that, upon mod $p$ reduction and group-change via $\sigma_1$,  $\RR(\Z)[1]$ specialises to $\RR(\F_p)[1]=\Z/p^2[1]$. This holds regardless of the value of $r$, because $\frob: \F_p \to \F_p$ is the identity. The rest of the proof is the same as that of Proposition \ref{proplifthom}, keeping track of the extra cyclotomic twist $.[1]$. Recall that, in this proof, $.(1)$  denotes the (geometric) twist by the line bundle $\mathcal O_{\P(V)}(1)$. The cyclotomic and geometric twists, are of very different origins.
\end{dem}

\bibliographystyle{plain}
\bibliography{biblitex.bib}

\end{document}